\documentclass[times,sort&compress,3p]{elsarticle}
\journal{Journal of Multivariate Analysis}
\usepackage[labelfont=bf]{caption}

\usepackage{amsmath,amsfonts,amssymb,amsthm,booktabs,color,epsfig,graphicx,hyperref,url}

\usepackage{color}
\usepackage{multirow}

\usepackage{hyperref}
\usepackage[latin1]{inputenc}
\usepackage{bbm}
\usepackage{dsfont}
\usepackage{enumerate,comment}
\usepackage{xspace}
\usepackage{latexsym}
\usepackage{amsfonts}
\usepackage{amsmath}
\usepackage{amssymb}
\usepackage{amsxtra}
\usepackage{amsthm}
\usepackage{url}
\usepackage{rotating}
\usepackage[bottom]{footmisc}
\usepackage{mathdots}
\usepackage{nameref}
\usepackage{hhline}
\usepackage{verbatim}
\usepackage{lmodern}
\theoremstyle{plain}                    % default
\newtheorem{thm}{Theorem}
\newtheorem{lem}{Lemma}

\newtheorem{prop}{Proposition}

\theoremstyle{definition}

\newtheorem{remark}{Remark}

\theoremstyle{remark}
\def\N{{\mathbb N}}
\def\R{{\mathbb R}}
\def\C{{\mathbb C}}
\def\Z{{\mathbb Z}}

\def\ninfty{\mathop{\longrightarrow}\limits_{T\to\infty}}

\newcommand{\bd}{\begin{displaymath}}
\newcommand{\ed}{\end{displaymath}}
\newcommand{\bea}{\begin{eqnarray}}
\newcommand{\eea}{\end{eqnarray}}
\newcommand{\bean}{\begin{eqnarray*}}
\newcommand{\eean}{\end{eqnarray*}}

\newcommand\be{\begin{eqnarray}}
\newcommand\ee{\end{eqnarray}}
\newcommand\bee{\begin{eqnarray*}}
\newcommand\eee{\end{eqnarray*}}

\begin{document}

\begin{frontmatter}

\title{Testing the equality of spectral density operators for functional  processes}

\author[A1]{Anne Leucht}
\author[A2]{Efstathios Paparoditis}
\author[A3]{Daniel Rademacher}
\author[A4]{Theofanis Sapatinas\corref{mycorrespondingauthor}}

\address[A1]{University of Bamberg, Research Group of Statistics and Business Mathematics, Feldkirchenstra\ss e 21, D-96052 Bamberg, Germany}
\address[A2]{University of Cyprus, Department of Mathematics and Statistics, P.O.Box 20537, CY-1678, Nicosia, Cyprus}
\address[A3]{Technische Universit\"at Braunschweig, Institut f\"ur Mathematische Stochastik, Universit\"atsplatz 2, D-38106 Braunschweig, Germany}
\address[A4]{University of Cyprus, Department of Mathematics and Statistics, P.O.Box 20537, CY-1678, Nicosia, Cyprus}

\cortext[mycorrespondingauthor]{Corresponding author. Email address: fanis@ucy.ac.cy \url{}}

\begin{abstract}
The problem of comparing the entire second order structure  of two functional processes is considered and 
a  $L^2$-type statistic for testing equality of the corresponding spectral density operators  is investigated. 
The test statistic  evaluates, over all frequencies, the Hilbert-Schmidt  distance 
between  the  two estimated spectral density operators. Under  certain assumptions,  the   limiting 
 distribution  under the null 
hypothesis  is derived.  A novel frequency domain bootstrap method is introduced,  
which  leads to a  more accurate approximation of the  distribution of  the test statistic under the null than 
the large sample Gaussian 
approximation derived.  Under quite general conditions, 
asymptotic validity  of the bootstrap procedure is established for  estimating  the distribution of the test statistic under the null. Furthermore, consistency of the bootstrap-based test under the alternative 
is  proved. Numerical simulations show that, even for small samples, 
 the bootstrap-based test has a very good size and power  behavior. 
An application to a  bivariate real-life  functional time series illustrates the methodology proposed.   
\end{abstract}


\end{frontmatter}

\section{Introduction}
\label{S1}
%%%%%%%%%%%%%%%%%%%%%%%%%%%%%%%%%%%%%%%%%%%%%%%%%%%%%%%%%%%%%%%%%%%%%%%%%%%%%%%

Functional data analysis  is a branch of statistics that in recent years has grown considerably and has created great research interest in the scientific community, especially in connection with the increasing number of situations in which theoretical and applied scientists have to deal with data of a continuous nature (i.e., curves, images, surfaces, etc.). For various works and references in different branches of functional data analysis, we refer to the recent special issues of Goia and Vieu \cite{GV-2016} and Aneiros {\rm et al.} \cite{AFV-2019}. See also the monograph by Horv\'ath \& Kokoszka \cite{HK-2012} which discusses inference problems in a variety of setting concerning independent as well as dependent functional data.

%\vspace*{0.3cm}
In our work, we focus on dependent functional data and, in particular, on functional time series analysis. Functional time series occurs  in many applications such as daily curves of financial transactions, daily images of geophysical and environmental data and daily curves of temperature measurements. Such curves or images are viewed as functions in appropriate spaces since an observed intensity is available at each point on a line segment, a portion of a plane or a volume. Moreover, and most importantly, such functional time series  exhibit temporal dependence and ignoring this  dependence  may result in misleading conclusions and not approperiate inferential procedures.
%\vspace*{0.3cm}

Comparing characteristics of two or more groups of functional data forms an important problem of statistical inference with a variety of applications. For instance, comparing the mean functions between independent groups of independent and identically distributed (i.i.d.) functional data  has attracted considerable interest in the  literature, see, e.g., Benko {\rm et al.} \cite{BHK-2009}, Zhang {\rm et al.} \cite{ZPZ-2010}, Horv\'ath and Kokoszka \cite{HK-2012} (Chapter 5), Horv\'ath {\rm et al.} \cite{HKR-2013} and Paparoditis and Sapatinas \cite{PS-2016}. In contrast to comparing  mean functions,  the problem of comparing the entire second order structure of two independent functional time series has been much less investigated.  Notice that for i.i.d. functional data this problem simplifies  to the problem of testing the equality of (the lag zero) covariance operators, see, e.g.,  Panaretos {\rm et al.} \cite{PKM-2010}, Fremdt {\rm et al.} \cite{FSHK-2012}, Pigoli {\rm et al.} \cite{PADS-2014} and Paparoditis and Sapatinas \cite{PS-2016}. The same problem of testing the equality of the (lag-zero) covariance operators of two sets of independent functional time series has also been investigated by Zhang and Shao \cite{ZS-2015}  and by Pilvakis {\rm et al.} \cite{PPS-2020}.
%\vspace*{0.3cm}

However, the comparison of the entire second order structure of independent functional time series,  is a much more involved problem due to the temporal dependence between 
  the  random elements considered.  In describing the second order structure of functional time series, the spectral density operator, introduced  in the functional set-up by Panaretos and Tavakoli \cite{PT-2013}, is  a very useful tool since it  summarizes in a nice way  the entire autocovariance structure of the underlying functional time series; see also H\"ormann {\rm et al.} \cite{HKH-2015}.
  % and van Delft and Eichler (2018).  
  It is, therefore,  very appealing  to develop a  spectral approach for testing  equality of the 
  entire second  order structure  of two functional time series. Tavakoli and Panaretos \cite{TP-2016} proposed an approach based on   projections on finite dimensional spaces of the differences of the estimated  spectral density  operators of the two functional time series. Projection-based tests  have the advantage to lead to manageable limiting distributions. However,  such tests  have no power  for  alternatives which are orthogonal to the projection space considered. Furthermore,  
  the number of projections appears as an additional tuning parameter  which has to be chosen by the user. Finally,  simulations in the much simpler i.i.d.~set-up suggest  that the quality of the large sample Gaussian  approximations of the corresponding test,  is affected by  the number of projections used; see Paparoditis and Sapatinas \cite{PS-2016}.  In this paper we focus  on tests which evaluate the differences between 
the entire, infinite dimensional, structure of the two spectral density operators compared.  For this, the Hilbert-Schmidt norm of the differences between  the estimated spectral density operators,  evaluated over all frequencies,  is used as the  basic building block of the test statistic considered. 
%Similar tests in the finite dimensional case have been proposed Notice that in the finite-dimensional case, i.e., for (univariate or multivariate) real-valued time series, such tests %have been developed among others by Eichler (2008) and  Dette and Paparoditis (2009).

The contribution of this  paper is twofold.   First, we   focus on  testing the equality of the entire second order structure between two independent  functional processes
by evaluating  for each frequency, the Hilbert-Schmidt norm between the (estimated) spectral density operators of the  functional  process at hand. Integrating these differences over all possible frequencies, leads to a global, $L^2$-type,  measure of deviation which is used to test the null hypothesis of interest. We show that 
under the assumption of a linear Hilbertian processes,  the limiting distribution of an appropriately  centered version of such a  test statistic under the null, is Gaussian.  This Gaussian distribution does not depend on characteristics of the underlying functional  processes  beyond those of second order. Second, and because of the slow convergence of the distribution of  the considered $L^2$-type test statistic  under the null against the derived limiting Gaussian distribution, we develop a novel frequency domain bootstrap procedure to estimate this distribution. The frequency domain bootstrap method works under minimal  conditions on the underlying functional process and its range of applicability is not restricted to the particular class of processes 
considered  and which is used  to derive the limiting distribution of our test. We prove under very  general conditions, 
 that the bootstrap procedure correctly approximates  the distribution of the proposed test statistic under the null. Furthermore, consistency of the bootstrap-based test under the alternative is established. Our theoretical deviations are accompanied  by a simulation study which shows a very good behavior of the bootstrap procedure in approximating the distribution of interest and  the good size and power performance of the test based on bootstrap critical values.  Notice that the frequency domain bootstrap method proposed  in this paper, can potentially be used to improve the performance of other tests too, like for instance, the projections based test of  Tavakoli and Panaretos \cite{TP-2016}. 

Developing bootstrap procedures for functional time series has attracted considerable  interest in the literature.  Politis and Romano \cite{PR-1994} established weak convergence 
results for the stationary  bootstrap,  Dehling {\rm et al.} \cite{DSW-2015} for the (non-overlapping) block bootstrap in a testing context, Ra\~na {\rm et al.} \cite{RAV-2015} applied a stationary bootstrap to functional time series,   Ferraty and Vieu \cite{FV-2011} a residual-based bootstrap and  Franke and Nyarige \cite{FN-2019} established consistency  of  a model-based bootstrap for  functional autoregressions. Pilavakis {\rm et al.} \cite{PPS-2019} derived  theoretical results for  the  moving block bootstrap and for the tapered  block bootstrap, Shang \cite{S-2018} applied  a maximum entropy bootstrap and Paparoditis \cite{P-2018} introduced a sieve  bootstrap for functional time series.  In contrast to the aforementioned contributions, the bootstrap procedure proposed in this paper acts  solely in the frequency domain and generates  replicates of the periodogram kernels stemming from functional processes that satisfy the null hypothesis of interest.

A  test related to ours and proposed after the first preprint of this paper has been appeared (see Leucht {\rm et al.} \cite{LPS-2018}), is that of  van Delft and Dette \cite{vDD-2020}, which  deals with  testing a different set of hypotheses, so-called relevant hypotheses, about the second order dynamics of two functional processes.  Important differences between the two procedures appear which will be discussed in more detail later on. 
However, we  stress here the fact that the test statistic proposed in this paper is not a special case of the test statistic used in the aforementioned paper  and, consequently,  the limiting distribution of our test statistic is different and not covered by the asymptotic  results derived in that paper. See  Remark~{\rm et al.} for more details.

The remainder  of the paper is organized as follows. Section 2 contains the main assumptions on the underlying functional linear processes and states the hypothesis testing problem under study. Section 3 is devoted to the suggested test statistic and its asymptotic behavior while Section 4 presents the frequency domain bootstrap procedure proposed to estimate the distribution of the test statistic under the null. Asymptotic validity of the bootstrap procedure is established and consistency of the  corresponding  test under the alternative also is proved. Section 5 contains numerical simulations and an application to  a bivariate meteorological functional time series while Section 6 concludes our findings. Auxiliary results 
containing some new results on frequency domain properties of linear Hilbertian processes as well as proofs of the main results are deferred to the Appendix and to the Supplementary Material.

%%%%%%%%%%%%%%%%%%%%%%%%%%%%%%%%%%%%%%%%%%%%%%%%%%%%%%%%%%%%%%%%%%%%%%%%%%%%%%%
\section{Assumptions and  the Testing Problem}
\label{S2}
%%%%%%%%%%%%%%%%%%%%%%%%%%%%%%%%%%%%%%%%%%%%%%%%%%%%%%%%%%%%%%%%%%%%%%%%%%%%%%%

Suppose that observations $X_1,\dots, X_T$ and $Y_1,\dots, Y_T$ stem from functional  processes $(X_t)_{t\in\Z}$ and $(Y_t)_{t\in\Z}$, respectively,  satisfying the following assumption.

%\vspace*{0.3cm}
\noindent
{\it Assumption 1}\ : 
 $(X_t)_{t\in\Z}$ and $(Y_t)_{t\in\Z}$ are independent functional linear processes, given by
\begin{equation}\label{eq.lin-proc}
X_t=\sum_{j\in \Z} A_j (\varepsilon_{t-j})\quad\text{and}\quad Y_t=\sum_{j\in \Z} B_j (e_{t-j}),\; t\in \Z,
\end{equation}
with values in $L^2_\R([0,1], \mu)$, where $\mu$ denotes the Lebesgue measure. The innovation functions  $(\varepsilon_t)_{t \in \Z}$ and $(e_t)_{t \in \Z}$ are two i.i.d.~mean zero Gaussian processes with values in $L^2_\R([0,1], \mu)$ and covariance operators $C_\varepsilon$ and $C_e$ with continuous
%\footnote{\textcolor{red}{DR: Is continouity of the kernels really needed in order to proof theorem 3.1?}} 
covariance kernels $ c_\varepsilon$ and $ c_e$, respectively. The sequences $(A_j)_{j \in \Z}$ and $(B_j)_{j \in \Z}$ of bounded linear operators from $L^2_\R([0,1], \mu)$ to $L^2_\R([0,1], \mu)$ where  $ A_0=B_0$ is the identity 
operator,
%\footnote{\textcolor{red}{DR: Is this actually necessary?}}, 
satisfy $\sum_{j \in \Z}|j|(\|A_j\|_{\mathcal L}+\|B_j\|_{\mathcal L})<\infty$ with $\|\cdot\|_{\mathcal L}$ denoting the operator norm.
%\vspace*{0.3cm}

%In particular, note that the summability conditions on the coefficients assure a.s.~absolute convergence of both series in \eqref{eq.lin-proc}.
We are interested in testing for equality of the entire second order structure of  the two functional  processes given in (\ref{eq.lin-proc}).  Notice that considering linear processes   in Assumption~1 should not be  considered as restrictive since  our interest is solely focused  on the  comparison of  the second order structure, i.e.,  of the  autocovariance  structure  of the  underlying functional processes. Furthermore,  and as we will see  later on, 
the  assumption of Gaussian innovation functions $ \varepsilon_t$ and $e_t $  is not essential. It  is solely imposed in order to  simplify  the already quite involved technical arguments 
used  to derive the   limiting  distribution of the test.  

For the  testing problem considered  it turns out that a spectral approach  is very appealing.
Towards this notice first that   we can define a spectral density operator in the sense of Panaretos and Tavakoli \cite{PT-2013} in the present set up which generalizes the concept of spectral densities for univariate time series and spectral density matrices for multivariate time series. Here and in the sequel, we will abbreviate $L^2_\R([0,1]^d,\mu)$ by $L^2$ if the dimension~$d$ becomes clear from the context.

\begin{lem}\label{l.ex-SDO}
Suppose that $(X_t)_{t\in\Z}$ and $(Y_t)_{t\in \Z}$ are functional processes satisfying Assumption~1. 
Then, for arbitrary $\lambda\in (-\pi,\pi]$,
$$
f_{X,\lambda}(\cdot,\cdot)=\frac{1}{2\pi}\sum_{t\in\Z}e^{-i\lambda t} \, r_{X,t}(\cdot,\cdot), \quad f_{Y,\lambda}(\cdot,\cdot)=\frac{1}{2\pi}\sum_{t\in\Z} {e^{-i\lambda t}}\, r_{Y,t}(\cdot,\cdot)
$$
with $r_{X,t}$ and $r_{Y,t}$ denoting the autocovariance kernels of $X$ and $Y$ at lag $t$, respectively, converge absolutely in $L^2$.
Moreover, for all $\sigma,\tau \in [0,1]$,
$$
r_{X,t}(\sigma, \tau)=\int_{(-\pi,\pi]} f_{X,\lambda}(\sigma,\tau)\,{e^{i\lambda t}}\, d\lambda, \quad
r_{Y,t}(\sigma, \tau)=\int_{(-\pi,\pi]} f_{Y,\lambda} (\sigma, \tau)\, {e^{i\lambda t}}\, d\lambda\quad
\forall t \in \mathbb{Z},
$$
where equality holds in $L^2$.  The operators ${\mathcal F_{X,\lambda}}$ and ${\mathcal F_{Y,\lambda}}$, induced by right integration of $f_{X,\lambda}$ and $f_{Y,\lambda}$, are self-adjoint, nonnegative definite and it holds 
$$
\mathcal F_{X,\lambda}=\frac{1}{2\pi}\sum_{t\in\Z} e^{-i\lambda t} \mathcal R_{X,t}, \quad\mathcal F_{Y,\lambda}=\frac{1}{2\pi}\sum_{t\in\Z} e^{-i\lambda t} \mathcal R_{Y,t},
$$ 
where $\mathcal R_{X,t}$ and $\mathcal R_{Y,t}$ denote the autocovariance operators of $X$ and $Y$ at lag $t$, induced by right integration of  $r_{X,t}$ and $r_{Y,t}$, respectively. Convergence holds in nuclear norm.
\end{lem}

%As it can be seen from the proof of Lemma \ref{l.ex-SDO}, it suffices to assume 4th order integrability of the innovations $(\varepsilon_t)_{t \in \Z}$ and $(e_t)_{t \in \Z}$ rather than %Gaussianity as in Assumption 1. 
The kernels $f_{X,\lambda}$ and $f_{Y,\lambda}$ are called the \textsl{spectral density kernels} (at frequency $\lambda$) and the operators $\mathcal F_{X,\lambda}$ and $\mathcal F_{Y, \lambda}$ are referred to as the corresponding {\sl spectral density operators}. 
%\vspace*{0.3cm}

Under the assumptions of Lemma~\ref{l.ex-SDO}, we can now state the hypothesis testing problem of interest as follows
\begin{equation}
\begin{aligned}
&{\mathcal H_0}\colon \mathcal F_{X,\lambda}=\mathcal F_{Y, \lambda} \quad\text{ for }\mu\text{-}\text{almost all } \lambda\in (-\pi,\pi], 
 \\
&{\mathcal H_1}\colon \mathcal F_{X, \lambda}\neq\mathcal F_{Y, \lambda} \quad \forall\lambda\in A \text{ for some $A\subset [0,\pi]$ with }\mu(A)>0.
\end{aligned}
\label{eq:testM}
\end{equation}

%%%%%%%%%%%%%%%%%%%%%%%%%%%%%%%%%%%%%%%%%%%%%%%%%%%%%%%
\section{The Test Statistic and its Asymptotic Behavior}\label{S3}
%%%%%%%%%%%%%%%%%%%%%%%%%%%%%%%%%%%%%%%%%%%%%%%%%%%%%%%

We first estimate the unknown spectral density operator $\mathcal F_{X, \lambda}$ by an integral operator~$\widehat{\mathcal F}_{X,\lambda}$ induced by right integration with the kernel
$$
\hat f_{X,\lambda}(\sigma,\tau)=\frac{1}{b T}\sum_{t=-N}^{N} W\left(\frac{\lambda-\lambda_t}{b}\right)\, \widehat p_{X, \lambda_t} (\sigma, \tau), \;\; \mbox{for all} \;\; \sigma,\tau \in [0,1],
$$
and, similarly, $\mathcal F_{Y, \lambda}$ by an integral operator~$\widehat{\mathcal F}_{Y,\lambda}$ induced by right integration with the kernel
$$
\hat f_{Y,\lambda}(\sigma,\tau)=\frac{1}{b T}\sum_{t=-N}^{N} W\left(\frac{\lambda-\lambda_t}{b}\right)\, \widehat p_{Y, \lambda_t} (\sigma, \tau),  \;\; \mbox{for all} \;\; \sigma,\tau \in [0,1].
$$
Here, $N=[(T-1)/2]$ and $\lambda_t=2\pi t/T, \;t \in \{-N,\dots, N\},$ denote the Fourier frequencies. Furthermore, $b=b_T>0$ is an asymptotically vanishing bandwidth and $W$ denotes a weight function. Moreover, as in Panaretos and Tavakoli \cite{PT-2013},
$$
\widehat p_{X, \lambda}(\sigma,\tau)=\frac{1}{2\pi T}\sum_{s_1,s_2=1}^T X_{s_1}(\sigma)X_{s_2}(\tau)\, \exp(-i\lambda( s_1- s_2)),  \;\; \mbox{for all} \;\; \sigma,\tau \in [0,1],
$$
and
$$
\widehat p_{Y, \lambda}(\sigma,\tau)=\frac{1}{2\pi T}\sum_{s_1,s_2=1}^T Y_{s_1}(\sigma)Y_{s_2}(\tau)\, \exp(-i\lambda( s_1- s_2)),  \;\; \mbox{for all} \;\; \sigma,\tau \in [0,1],
$$
denote the periodogram kernels based on $X_1,\dots, X_T$ and  $Y_1,\dots, Y_T$, respectively.  The periodogram  operators $I_{X,\lambda}$, and $I_{Y,\lambda}$ are defined as  integral operators  induced by right integration of the periodogram kernels $\widehat p_{X, \lambda}$ and $\widehat p_{Y, \lambda}$, respectively.  
%\vspace*{0.3cm}

For the hypothesis testing problem (\ref{eq:testM}), we propose the following test statistic 
\begin{equation}
{\mathcal U_T}=\int_{-\pi}^\pi\|\widehat{\mathcal F}_{X,\lambda}-\widehat{\mathcal F}_{Y,\lambda}\|_{HS}^2\, d\lambda,
\label{eq:test-stat-main}
\end{equation}
which evaluates the distance between the estimated spectral density operators via the Hilbert-Schmidt norm $ \|\cdot\|_{HS}$. 
%where $\|\cdot\|_{HS}$ denotes the Hilbert Schmidt norm on $L^2$. 
The following theorem states  the  asymptotic properties of the suitably normalized test statistic $\mathcal U_T$ when the null hypothesis $ {\mathcal H_0}$ is true.

%\vspace*{0.3cm}

\begin{thm}\label{t.asy-test}
Suppose that the stretches of observations $X_1,\dots, X_T$ and $Y_1,\dots, Y_T$ stem from the two  functional processes $(X_t)_{t\in\Z}$ and $(Y_t)_{t\in \Z}$, respectively,  satisfying Assumption~1. Moreover, assume that \\
%\begin{itemize}
%\item[] 
{\rm (i)} $b\sim T^{-\nu }$ for some $\nu\in (1/4,1/2)$, \\
%\item[(ii)] 
{\rm (ii)} $W$ is bounded, symmetric, positive, and Lipschitz continuous, has  bounded support on $(-\pi,\pi]$ and satisfies $\int_{-\pi}^{\pi} W(x)\, dx=2\pi$. \\
%\item[(iv)] $\int_0^1 E\varepsilon_0^8(\sigma)\, d\sigma<\infty$.
%\end{itemize}
Then, under ${\mathcal H_0}$,
\begin{equation} \label{eq.DistrH0}
\sqrt{b} T \, \mathcal U_T-b^{-1/2} \mu_0\stackrel{d}{\longrightarrow} Z\sim{\mathcal N}(0,\theta_0^2),
\end{equation}
where
$$
\begin{aligned}
\mu_0&=\frac{1}{\pi}\, \int_{-\pi}^\pi\left\{\mathrm{trace}(\mathcal F_{X,\lambda})\,\right\}^2 d\lambda\int_{-\pi}^\pi W^2\left(u\right)\, du,\\
\theta_0^2&=
%\textcolor[rgb]{1,0,0}{
\frac{4}{\pi^2}\int_{-2\pi}^{2\pi}\left\{\int_{-\pi}^\pi W(u)W(u-x)\, du\right\}^2dx\, 
\int_{-\pi}^\pi \| \mathcal F_{X,\lambda} \|_{HS}^4 \,d\lambda \, .	
%}
\end{aligned}
$$
\end{thm}

%\vspace*{0.3cm}

Note that the assumptions (i) and (ii) on the weight function $W$ and the bandwidths $(b_T)_T$, respectively, in Theorem \ref{t.asy-test} are identical to the assumptions for multivariate time series used in Dette and Paparoditis \cite{DP-2009}.

\begin{remark}
\label{rem:3.1}
{\rm In our work, we have considered the case where the sample sizes of both time series $(X_t)_{t=1}^T $ and  $(Y_t)_{t=1}^T$ are equal. In principle, we could also consider time series of different length, that is 
$(X_t)_{t=1}^{T_1}$ and $(Y_t)_{t=1}^{T_2}$. Under certain regularity conditions, such as $\sqrt{b_1}T_1/(\sqrt{b_1}T_1+\sqrt{b_2}T _2)\to \eta\in(0,1)$ as $T_1+T_2\to\infty$,  and with minor, but tedious modifications of the proof, one can also show asymptotic normality of  $(\sqrt{b_1}T_1+\sqrt{b_2}T _2) \, {\mathcal U}_{T_1,T_2}$, after a suitable centering. Here, ${\mathcal U}_{T_1,T_2}=\int_{-\pi}^\pi \|\widehat {\mathcal F}_{X,\lambda}^{(T_1)}-\widehat {\mathcal F}_{Y,\lambda}^{(T_2)}\|_{HS}^2\,d\lambda$  relies on the estimated spectral density operator $\widehat {\mathcal F}_{X,\lambda}^{(T_1)} $, based on $(X_t)_{t=1}^{T_1}$, and  the estimated spectral density operator $\widehat {\mathcal F}_{Y,\lambda}^{(T_2)} $, based on $(Y_t)_{t=1}^{T_2}$,  using bandwidths $b_1$ and $b_2$, respectively.}
\end{remark}

\begin{remark}
{\rm A careful  inspection of the proof of Theorem~\ref{t.asy-test} shows that the assumption of Gaussianity on the functional innovations  $(\varepsilon_t)_{t \in \Z}$ and $(e_t)_{t \in \Z}$  in   (\ref{eq.lin-proc}) is  solely used  to simplify somehow the technical arguments applied  in 
proving  asymptotic normality of the quadratic forms    involved in proving assertion (\ref{eq.DistrH0}) of Theorem~\ref{t.asy-test}.  Notice that  this   assumption is not required  in order to prove   convergence of the  mean and  of the variance of $ \sqrt{b}T \,\mathcal U_T$ to the limits given in the aforementioned theorem.  Consequently, this assumption  can be  replaced by     
other assumptions on the stochastic properties of  the innovations $(\varepsilon_t)_{t \in \Z}$ and $(e_t)_{t \in \Z}$,  which will allow  for the  use of different technical arguments, for instance arguments based on the convergence of  all cumulants of the random sequence $\sqrt{b} T \, \mathcal U_T-b^{-1/2} \mu_0$   to the appropriate limits, in order to establish the desired asymptotic normality. Furthermore, the bootstrap approach proposed in the next section does not rely on and it does not make use of  the structural assumptions imposed on the underlying functional processes in order to derive the limiting distribution of the test. }
\end{remark}

\begin{remark} \label{re.DeDe1}
A closely related null hypothesis $\mathcal H_0\colon \int_a^b\| {\mathcal F}_{X,\lambda}- {\mathcal F}_{Y,\lambda}\|_{HS}^2\,d\lambda \leq \Delta$ has been considered in van Delft and Dette \cite{vDD-2020} for prespecified constants $a<b\in[0,\pi]$ and $\Delta >0$. Although their test statistic proposed looks at a first glance similar  to ours, see equation (3.19) in the aforecited paper, several differences appear. Notice first that the convergence rate of the nominator and of the denominator of their statistic is of  order $ O_P(\sqrt{Tb})$ and not  $ O_P(T\sqrt{b})$,  as    of  the test statistic  (\ref{eq:test-stat-main}) considered in  this paper. Apart from the fact that a different set of null hypotheses is considered in the two papers, the main reason for this  difference in the convergence rates,  lies in the fact that the limiting distribution of the test statistic considered in van Delft and Dette \cite{vDD-2020} is essentially dominated 
by the differences $ \widehat{\mathcal F}_{X,\lambda} -{\mathcal F}_{X,\lambda}$, respectively,    $ \widehat{\mathcal F}_{Y,\lambda} -{\mathcal F}_{Y,\lambda}$, which are of order  $ \sqrt{Tb}$. On the other hand, the distribution of our  test statistic is dominated by the quadratic term  $ \| \widehat{\mathcal F}_{X,\lambda}-\widehat{\mathcal F}_{Y,\lambda}\|_{HS}^2 $, which   in the test statistic considered by van Delft and Dette \cite{vDD-2020} disappears; see Lemma 3.1 of their paper. Consequently, to establish asymptotic normality of the test statistic considered in van Delft and Dette \cite{vDD-2020}, essentially, a central limit theorem for $ \sqrt{Tb}(\widehat{\mathcal F}_{X,\lambda} -{\mathcal F}_{X,\lambda})$ , respectively, for $  \sqrt{Tb}(\widehat{\mathcal F}_{Y,\lambda} -{\mathcal F}_{Y,\lambda})$ is  involved. In contrast to this, our test statistic deals with  weighted sums of 
the quadratic  terms  $\Big< \widehat{\mathcal F}_{X,\lambda_1}-\widehat{\mathcal F}_{Y,\lambda_1}, \widehat{\mathcal F}_{X,\lambda_2}-\widehat{\mathcal F}_{Y,\lambda_2}\Big>_{HS} $,  for which central limit theorems for generalized quadratic forms has to be invoked.  Even in the finite dimensional case, central limit theorems  for generalized quadratic forms are established  under more structural  assumptions on the underlying processes than those needed to deal with  the sequence  $\sqrt{Tb}(\widehat{\mathcal F}_{X,\lambda} -{\mathcal F}_{X,\lambda})$; see for instance Eichler \cite{E-2008} who uses summability conditions on the cumulants of all order or Dette and Paparoditis \cite{DP-2009} who use linearity assumptions on the underlying vector processes. The technical challenges  in dealing with the test statistic (\ref{eq:test-stat-main}), also justify the additional structural assumptions imposed  in this  paper in order   to establish the  limiting distribution of $ {\mathcal U}_T$, as compared to those used  in van Delft and Dette \cite{vDD-2020}.
\end{remark}

%\vspace*{0.3cm}

Based on Theorem~\ref{t.asy-test}, the procedure to test hypothesis (\ref{eq:testM}) is then defined as follows:  Reject  ${\mathcal H_0}$ if and only if 
\begin{equation}    \label{eq.stud-test}
t_{\mathcal U} = \frac{\sqrt{b} T \, \mathcal U_T-b^{-1/2} \widehat{\mu}_0}{\widehat{\theta}_0} \geq z_{1-\alpha},
\end{equation}
where  $ z_{1-\alpha}$ is the upper $1-\alpha$ percentage point of the standard  Gaussian  distribution
and $ \widehat{\mu}_0$ and $\widehat{\theta}_0$ are consistent estimators of $ \mu_0$ and $\theta_0$, respectively. Such estimators 
can be, for instance,  obtained if  the unknown spectral density kernel $ f_{X,\lambda}$ is replaced by the  pooled estimator 
$ \widehat{f}_{\lambda}(\tau,\sigma) =\widehat{f}_{X,\lambda}(\tau,\sigma)/2 +  \widehat{f}_{Y,\lambda}(\tau,\sigma)/2$.
%given in (\ref{eq.sp-pooled}) in Section~\ref{sec.boot}.  
Notice that,  under ${\mathcal H_0}$, $ f_{X,\lambda}=f_{Y,\lambda}=f_{X,\lambda} /2+f_{Y,\lambda} /2$, that is (asymptotically), it makes no difference  if $ f_{X,\lambda}$ in $\mu_0$ and $\theta_0$ is replaced by $ \widehat{f}_{X,\lambda}$ (or by $ \widehat{f}_{Y,\lambda}$) instead of the pooled estimator $\widehat{f}_\lambda$. However, under ${\mathcal H_1}$ it matters and, for this reason, we  use the pooled estimator $\widehat{f}_\lambda(\tau,\sigma)$ in applying the studentized test statistic $t_{\mathcal U}$ defined in (\ref{eq.stud-test}); see also Lemma~\ref{le.power} in Section 4.  
Under the assumption that the pooled estimator $\widehat{f}_\lambda$ is uniformly consistent,  (see also Assumption 2 below),
%\[ \sup_{(\tau,\sigma) \in [0,1]^2}\sup_{\lambda\in [0,\pi]} \big|\widehat{f}_{\lambda}(\tau,\sigma) - f_{\lambda}(\tau,\sigma) \big|  = o_P(\sqrt{b}),\]
%as $ T \rightarrow \infty$, where $ f_\lambda(\tau,\sigma) =f_{X,\lambda}(\tau,\sigma)/2 +f_{Y,\lambda}(\tau,\sigma)/2$ for all $ (\tau,\sigma) \in [0,1]^2$. 
it is  easily seen that, under ${\mathcal H_0}$, 
\[ t_{\mathcal U}  =\frac{\sqrt{b} T \, \mathcal U_T-b^{-1/2} \mu_0}{\theta_0} +o_P(1), \]
i.e., Theorem~\ref{t.asy-test}  implies that the studentized test $t_{\mathcal U} $ is  an asymptotically $\alpha$-level test under ${\mathcal H_0}$, for any desired  level $\alpha \in (0,1)$.

\begin{remark}
Notice that the test statistic  $t_{\mathcal U}$ is asymptotically pivotal, i.e., its distribution under the null does not depend on any unknown characteristics of the underlying functional processes. Furthermore, the denumerator $ \theta_0$   can be estimated using the estimators of the spectral density operators involved in calculating the test statistic  ${\mathcal U}_T$. A problem, however,   occurs from the well-known fact that, even in the finite-dimensional case, the convergence of the distribution of such $L^2$-norm based tests towards their limiting  (Gaussian) distribution is very slow;  see, e.g.,  H\"ardle and Mammen \cite{HM-1993}, Paparoditis \cite{P-2000} and Dette and Paparoditis \cite{DP-2009}. In this case, bootstrap-based  approaches may be very effective. This issue is addressed in the next section where a frequency domain bootstrap procedure is developed 
and  its asymptotic validity is established.
\end{remark}
%\underbrace{\sum^m_{j=1}\xi_{j,t}v_j}_{\displaystyle X_{t,m}}
\section{Bootstrapping The Test Statistic} \label{sec.boot}

%A problem in implementing the above test occurs from the well-known fact that, even in the finite-dimensional case, the convergence of such $L^2$-norm based test statistics towards their limiting  distribution is very slow;  see, e.g.,  H\"ardle and Mammen (1993), Paparoditis (2000) and Dette and Paparoditis (2009). In this case, bootstrap-based  approaches may be very effective.  

In this section  we  propose a  novel frequency domain bootstrap procedure which can be used to estimate  the distribution of the test statistic $\mathcal U_T$ defined in (\ref{eq:test-stat-main}) and, of the studentized test $ t_{\mathcal U}$ defined in (\ref{eq.stud-test}) under ${\mathcal H_0}$.    The frequency domain bootstrap approach proposed is of interest on its own and can potentially be applied to other  test statistics or  testing problems  developed for  comparing frequency domain characteristics of the  functional processes.

We begin by recalling the fact that  for any $k\in \N$ and any   set of  points $ 0\leq s_1 < s_2 < \cdots < s_k\leq 1$ in the interval $[0,1]$, 
the corresponding $k$-dimensional  vector 
 of finite Fourier transforms 
 $$ J_{X,\lambda}=\Big(J_{X,\lambda}(s_j) = (2\pi T)^{-1/2}\sum_{t=1}^{T} X_t(s_j) e^{-i t\lambda},\ j  \in \{1,2, \ldots, k\}\Big), $$ 
 satisfies  for $ \lambda \in (0,\pi)$, 
\begin{equation} \label{eq.fft1}
\left( \begin{array}{c} J_{X,\lambda}(s_1)\\J_{X,\lambda}(s_2)\\ \vdots\\ J_{X,\lambda}(s_k) \end{array}\right) 
\stackrel{d}{\rightarrow} \mathcal N_C\Big( \left( \begin{array}{c} 0\\0\\ \vdots \\ 0 \end{array}\right), \underbrace{\left( \begin{array}{cccc} f_{X,\lambda}(s_1,s_1) &    f_{X,\lambda}(s_1,s_2) &  \ldots &  f_{X,\lambda}(s_1,s_k) \\
 f_{X,\lambda}(s_2,s_1) &    f_{X,\lambda}(s_2,s_2) & \ldots &  f_{X,\lambda}(s_2,s_k) \\ 
  \vdots & \vdots & \ldots & \vdots \\ f_{X,\lambda}(s_k,s_1) &    f_{X,\lambda}(s_k,s_2)& \ldots &  f_{X,\lambda}(s_k,s_k)\end{array} \right)}_{\displaystyle =\Sigma_\lambda} \Big),
\end{equation}
where $\mathcal N_C$ denotes a circularly-symmetric complex Gaussian  distribution with mean zero and complex-valued covariance matrix $ \Sigma_\lambda$. Furthermore, for two different 
frequencies  $0 < \lambda_j \neq \lambda_k <\pi$, the corresponding vectors  of finite Fourier transforms $J_{X,\lambda_j}$ and $J_{X,\lambda_k} $ are asymptotically independent; see, e.g.,  Theorem 5 in Cerovecki and H\"ormann \cite{CH-2017}.  These properties of $ J_{X,\lambda} $ and $ J_{Y,\lambda} $  as well as the fact that
$ \widehat{p}_{X,\lambda} (\sigma,\tau) = J_{X,\lambda}(\sigma)\overline{J}_{X,\lambda}(\tau)$, for  $\sigma, \tau \in [0,1]$, is  the periodogram kernel,  motivate the following bootstrap procedure to approximate the distribution of the test statistic $\mathcal U_T$ defined in (\ref{eq:test-stat-main}) under $H_0$.

\vspace*{0.2cm}
\begin{enumerate}
 \item[{\rm Step 1}:] \ For ${\lambda_t=2\pi t/T}$, $t \in \{1, \ldots, N\}$,   $ N=[(T-1)/2]$, estimate the  pooled spectral density operator ${ \mathcal F_{\lambda_t}}$ by 
\begin{equation}
\label{eq.sp-pooled}
  \widehat{\mathcal F}_{\lambda_t} = \frac{1}{2}\widehat{\mathcal F}_{X,\lambda_t} +  \frac{1}{2}\widehat{\mathcal F}_{Y,\lambda_t}
 \end{equation}
and denote  by $ \widehat{f}_{\lambda_t}(\sigma,\tau)$, for $ \sigma, \tau \in \{s_1, \ldots, s_k\}$, the corresponding estimated pooled spectral density kernel.
\item[{\rm Step 2}:] \ Generate two independent vectors  $J^\ast_{X,\lambda_t} $ and $J^\ast_{Y,\lambda_t} $ as 
\[ J^\ast_{X,\lambda_t} \sim \mathcal N_C(0, \widehat{\Sigma}_{\lambda_t})  \ \ \mbox{and} \ \  J^\ast_{Y,\lambda_t} \sim \mathcal N_C(0,  \widehat{\Sigma}_{\lambda_t}),\]
independently for $\lambda_1,\dots, \lambda_N$, where $ \widehat{\Sigma}_\lambda $ is  the matrix  obtained by replacing in $ \Sigma_\lambda$ the unknown spectral density kernel $ f_{X,\lambda}$ by its pooled estimator $ \widehat{f}_\lambda$.
For $\sigma, \tau \in \{s_1,\ldots, s_k\}$, let 
\[ p^\ast_{X,\lambda_t}(\sigma,\tau) = J^\ast_{X,\lambda_t}(\sigma)\overline{J}^\ast_{X,\lambda_t}(\tau)  \ \ \mbox{and} \ \ p^\ast_{Y,\lambda_t}(\sigma,\tau) = J^\ast_{Y,\lambda_t}(\sigma)\overline{J}^\ast_{Y,\lambda_t}(\tau)\]
while, for  $ t \in \{-1, \ldots, -N\}$, set  
$$ p^\ast_{X,\lambda_t}(\sigma,\tau)=\overline{p}^\ast_{X,-\lambda_t}(\sigma, \tau) \;\; \mbox{and} \;\;  p^\ast_{Y,\lambda_t}(\sigma,\tau)=\overline{p}^\ast_{Y,-\lambda_t}(\sigma,\tau).
$$ Furthermore, set for simplicity $ J^\ast_{X,0} = J^\ast_{Y,0}=0$.
\item[{\rm Step 3}:] \  For $\sigma, \tau \in \{s_1,\ldots, s_k\}$, let 
\[ \hat f^\ast_{X,\lambda_t}(\sigma,\tau)=\frac{1}{b T}\sum_{s=-N}^{N} W\left(\frac{\lambda_t-\lambda_s}{b}\right)\, \widehat p^\ast_{X, \lambda_s} (\sigma, \tau)\] 
and
\[ \hat f^\ast_{Y,\lambda_t}(\sigma,\tau)=\frac{1}{b T}\sum_{s=-N}^{N} W\left(\frac{\lambda_t-\lambda_s}{b}\right)\, \widehat p^\ast_{Y, \lambda_s} (\sigma, \tau). \]
\item[{\rm Step 4}:]  \ 
%Denote by   $\widehat{F}^\ast_{X,\lambda}$ and $ \widehat{F}_{Y,\lambda}$ the integral operators induced by right-integration with kernels
%$\hat f^\ast_{X,\lambda}(\sigma,\tau) $ and $\hat f^\ast_{Y,\lambda}(\sigma,\tau) $ respectively. 
Approximate the distribution of the test statistic $ \mathcal U_T $ defined in (\ref{eq:test-stat-main}) by the distribution of the  bootstrap test statistic $ \mathcal U^\ast_{T,k} $ given by 
\[  \mathcal U^\ast_{T,k}=\frac{2\pi}{Tk^2}\sum_{\substack{l=-N }}^{N}\sum_{i,j=1}^k \Big| \widehat{f}^\ast_{X,\lambda_l}(s_i,s_j)
-\widehat{f}^\ast_{Y,\lambda_l}(s_i,s_j) \Big|^2.\]  
\end{enumerate}

\begin{remark}
\label{rem:4.1}
{\rm The set of points $ 0\leq s_1 < s_2 < \cdots < s_k\leq 1
$  at which the $k$-dimensional complex-valued random vectors  
  $ J^\ast_{X,\lambda_t}$ and $ J^\ast_{Y,\lambda_t}$ are generated    can be set equal to the 
  set of sampling points  at which the  functional  random elements  $ X_t $ and $ Y_t $ are observed in reality. 
However, and as it is commonly done in functional data analysis, 
these finite-dimensional vectors  can be  transformed to functional objects using 
%the generated random vectors $  J^\ast_{X,\lambda_t}$ and $ J^\ast_{Y,\lambda_t}$ and  
a basis in  $L^2$, for instance,  the  Fourier basis. In this case, the bootstrap approximation of the test statistic $ \mathcal U_T$ defined in (\ref{eq:test-stat-main}) will then be given by  
\begin{equation}
\mathcal U^\ast_T=\frac{2\pi}{T}\sum_{\substack{l=-N}}^{N}\int_0^1\int_0^1 \Big| \widehat{f}^\ast_{X,\lambda_l}(\tau,\sigma)
-\widehat{f}^\ast_{Y,\lambda_l}(\tau,\sigma) \Big|^2d\tau d\sigma =  \frac{2\pi}{T}\sum_{\substack{l=-N}}^{N}\|\widehat{\mathcal F}_{X,\lambda_l}^\ast-
\widehat{\mathcal F}_{Y,\lambda_l}^\ast\|^2_{HS}.
\label{eq:boot-test-stat}
\end{equation}
 From an asymptotic point of view
 %,  and as an inspection of the proof of Theorem~\ref{t.boo-test} below shows, 
 both bootstrap approximations, $ \mathcal U^\ast_{T,k}$ and $ \mathcal U^\ast_{T} $, will lead to the same result, provided that for $ \mathcal U^\ast_{T,k}$ the number of points $ k$ increases to infinity as the  sample size $T$ increases to infinity. In our theoretical derivations we will concentrate on $ \mathcal U^\ast_{T}$.} 
\end{remark}
%\vspace*{0.3cm}

\begin{remark}
{\rm In the case where the sample sizes of both time series $(X_t)_{t=1}^{T_1} $ and  $(Y_t)_{t=1}^{T_2}$ are different (see Remark \ref{rem:3.1}), the bootstrap algorithm can be adapted accordingly. In particular, the estimated pooled spectral density operator $\widehat{\mathcal F}_{\lambda}$, used in Step 1 above,  can be obtained for any frequency $\lambda \in [0, \pi]$ as 
$$
\widehat{\mathcal F}_{\lambda} = \frac{T_1}{T_1+T_2} \widehat {\mathcal F}_{X,\lambda}^{(T_1)} +  \frac{T_2}{T_1+T_2}\widehat {\mathcal F}_{Y,\lambda}^{(T_2)},
$$
where the estimated spectral density operators $\widehat {\mathcal F}_{X,\lambda}^{(T_1)} $ and $\widehat {\mathcal F}_{Y,\lambda}^{(T_2)} $ are given in Remark \ref{rem:3.1}. Then, $J^\ast_{X,\lambda_{t_1}} $ and $J^\ast_{Y,\lambda_{t_2}}$ can be generated as in Step 2, but for the Fourier frequencies $\lambda_{t_1}$ and $\lambda_{t_2}$ corresponding to the sample sizes $T_1$ and $T_2$, respectively. Although a bootstrap version of the test statistic ${\mathcal U}_{T_1,T_2}$ given in Remark ~\ref{rem:3.1} can be defined,
the theoretical derivations to establish bootstrap consistency in this case are more involved and beyond the scope of this paper. }
\end{remark}

Following the bootstrap procedure described in Steps 1-4, a bootstrap-based test then rejects ${\mathcal H_0}$ if 
\[ 
%\frac{\sqrt{b} T \, \mathcal U_T-b^{-1/2} \widehat{\mu}_0}{\widehat{\theta}_0} 
t_{\mathcal U} \geq t^\ast_{\mathcal U,1-\alpha},\]
where  $ t^\ast_{\mathcal U, 1-\alpha}$  denotes  the upper $1-\alpha$ percentage point of the distribution of the bootstrap studentized test 
\begin{equation}
\label{tst-fanis}
t^\ast_{\mathcal U}=(\sqrt{b} T\,\mathcal U^\ast_T  -b^{-1/2}\widehat{\mu}^\ast_0)/\widehat{\theta}^\ast_0,
\end{equation} 
where 
$\mathcal U^\ast_T$ is defined in (\ref{eq:boot-test-stat}) and $\widehat{\mu}_0^\ast$ and $ \widehat{\theta}^\ast_0$ are obtained by replacing the unknown spectral density kernel $ f_{X,\lambda}$  in the expressions for $\mu_0$ and $\theta_0$ given in Theorem \ref{t.asy-test} by its pooled estimator 
$ \widehat{f}^\ast_{\lambda}(\sigma, \tau) = \widehat{f}^\ast_{X,\lambda}(\sigma, \tau)/2+ \widehat{f}^\ast_{Y,\lambda}(\sigma, \tau)/2 $, for all $\sigma, \tau \in [0,1]$.
Notice that this distribution can be evaluated by Monte Carlo.

%\vspace*{0.3cm}
\begin{remark}
{\rm It is worth mentioning that,  by the definition of $ \widehat{\mu}^\ast_0$  and $ \widehat{\theta}^\ast_0$, the  bootstrap studentized test $t^\ast_{\mathcal U}$ imitates correctly also the randomness  in $t_{\mathcal U}$ which is  introduced by replacing the unknown spectral density kernel  $ f_{X,\lambda}$ appearing in $ \mu_0$ and $\theta_0$ by its pooled estimator $ \widehat{f}_\lambda$; see (\ref{eq.stud-test}).  A computationally simpler alternative will be to ignore this asymptotically negligible effect, that is, to 
use, instead of $t^\ast_{\mathcal U}$ given in (\ref{tst-fanis}), the studentized version   $t^+_{\mathcal U}=(\sqrt{b} T\,\mathcal U^\ast_T  -b^{-1/2}\widehat{\mu}_0)/\widehat{\theta}_0$ of the bootstrap-based test.}
\end{remark}
%\vspace*{0.3cm}

Before describing  the asymptotic behavior of the bootstrap test statistic  $ \mathcal U^\ast_T$ defined in (\ref{eq:boot-test-stat}),  we  state the following assumption which clarifies our requirements  on the pooled spectral density kernel estimator $ \widehat{f}_\lambda$ used.

%\vspace*{0.3cm}
{\it Assumption 2}\ : The pooled spectral density kernel estimator $ \widehat{f}_\lambda $ satisfies
%\[ \sup_{(\tau,\sigma) \in [0,1]^2}\sup_{\lambda\in [0,\pi]} \big|\widehat{f}_{\lambda}(\tau,\sigma) - f_{\lambda}(\tau,\sigma) \big|  = o_P(\sqrt{b}), \quad \mbox{as} \;\; T \rightarrow \infty, \]
 \[ \sup_{\lambda_t\in \{2\pi k/T\mid k=1,\dots, N\}} \Big|\int_0^1 \int_0^1 \big(\widehat{f}_{\lambda_t}(\sigma,\tau) - f_{\lambda_t}(\sigma,\tau)\big)d\sigma d\tau \Big|  = o_P(\sqrt{b}), \quad \mbox{as} \;\; T \rightarrow \infty, \]
where $ f_\lambda$ is the spectral density kernel of the pooled spectral density operator $\mathcal F_\lambda=(1/2) \mathcal F_{X,\lambda} + (1/2) \mathcal F_{Y,\lambda} $.

%\vspace{0.3cm}
Notice that the above assumption can be easily verified by using results for uniform consistency of spectral density estimators of univariate time series, since 
\begin{align*}
\int_0^1\int_0^1 \widehat{f}_{X,\lambda} (\sigma,\tau)d\sigma d\tau & = \frac{1}{Tb}\sum_{t=-N}^N W\Big(\frac{\lambda-\lambda_t}{b}\Big) \int_0^1 \int_0^1 \widehat{p}_{X,\lambda_t}(\sigma,\tau) d\sigma d\tau
\end{align*}
can be interpreted as a kernel estimator of the spectral density of the univariate time series $ \int_0^1 X_t(s)ds$, $ t=1, \ldots, n$,  the periodogram of which at frequency $ \lambda_t$ equals  $\int_0^1 \int_0^1 \widehat{p}_{X,\lambda_t}(\sigma,\tau) d\sigma d\tau$.  For instance, for the linear functional process $\{X_t, t\in\Z\}$ considered in this paper, $\int_0^1 X_t(s)ds$ is a univariate linear process as well and, under certain conditions,  Assumption 2 is satisfied; see  Franke and H\"ardle \cite{FH-1992}. {Assumption 2  can also be fulfilled} under  different conditions on the integrated process $\int_0^1 X_t(s)ds$; see Wu and Zaffaroni \cite{WZ-2015} for a discussion. \\

The following theorem establishes the asymptotic validity of the suggested bootstrap procedure.\\

\begin{thm}\label{t.boo-test}
Suppose that  Assumptions~2 as well as the conditions (i) and (ii) of  Theorem~\ref{t.asy-test} are satisfied. 
Then, conditional on $X_1,\ldots,X_T$ and $Y_1,\ldots, Y_T$, as $T\rightarrow \infty$, 
\[ \sqrt{b} T \, \mathcal U^\ast_T-b^{-1/2} \widetilde{\mu}_0  \stackrel{d}{\rightarrow}  \mathcal N(0,\widetilde{\theta}_0),\]
in probability, where 
\begin{align*}
 \widetilde{\mu}_0 &= \frac{1}{\pi}\, \int_{-\pi}^\pi
\{\mathrm{trace}({\mathcal F}_\lambda)\}^2 d\lambda\int_{-\pi}^\pi W^2\left(u\right)\, du,\\
\widetilde{\theta}_0^2 &=\frac{4}{\pi^2}\int_{-2\pi}^{2\pi}\left\{\int_{-\pi}^\pi W(u)W(u-x)\, du\right\}^2dx\, \int_{-\pi}^\pi 
\|{\mathcal F}_\lambda\|_{HS}^4\, d\lambda
\end{align*} 
and ${\mathcal F}_\lambda$ is the pooled spectral density operator given in  Assumption 2.
\end{thm}

%\begin{enumerate}
%\item[(i)] if  $ H_0$ is true,
%\[ d_K\big(\sqrt{b} T \, \mathcal T_T-b^{-1/2} \mu_0,\sqrt{b} T \, \mathcal T^\ast_T-b^{-1/2} \mu_0 \big) \rightarrow 0,\]
%in probability, while
%\item[(ii)]  if $H_1$ is true
%$$d_K\big( \big) \rightarrow 0,$$
%in probability, where 
%\[ \widetilde{\mu}_0 = \frac{1}{\pi}\, \int_{-\pi}^\pi\left\{\int_0^1 f_{\lambda}(\tau,\tau) \,d\tau\,\right\}^2 d\lambda\int_{-\pi}^\pi W^2\left(u\right)\, du,\]
%\begin{align*}
%\widetilde{\theta}_0^2&=\frac{2}{\pi\textcolor[rgb]{1,0,0}{^2}}\int_{-2\pi}^{2\pi}\left\{\int_{-\pi}^\pi W(u)W(u-x)\, du\right\}^2dx\, \int_{-\pi}^\pi \left\{\int_0^1\int_0^1f_{X,\lambda}^2(\sigma, \tau)\,\,d\sigma\, d\tau\right\}^2\, d\lambda
%\end{align*}
%and  $ f_\lambda $ is the kernel of the integral operator $(1/2) F_{X,\lambda} + (1/2) F_{Y,\lambda} $.
%\end{enumerate}

Notice that, under $\mathcal H_0$, $ \mu_0=\widetilde{\mu}_0$  and 
$\theta^2_0=\widetilde{\theta}_0^2$ since $ \mathcal F_{X, \lambda}= \mathcal F_{Y, \lambda} $ (or, respectively, $ f_{X,\lambda}=f_{Y, \lambda}$). Thus, in this case,  the asymptotic behavior of the test statistics $\mathcal U_T$ and $\mathcal U^\ast_T$ is identical, that is,  the bootstrap procedure estimates consistently the distribution of the test statistic $ \mathcal U_T$ under $\mathcal H_0$. Furthermore,  under $\mathcal  H_1$, the following holds true.

\begin{remark}
{\rm  As Theorem~\ref{t.boo-test}  shows,  the limiting distribution of the appropriately centered  bootstrap test statistic $ {\mathcal U}^\ast_T $ is obtained under validity of Assumption 2 and without imposing any particular assumptions on the weak dependence structure of the underlying functional processes $ \{ X_t, t \in \Z\}$ and  $ \{ Y_t, t \in \Z\}$. That is, this bootstrap procedure will lead to (asymptotically) valid approximations  for the same test if 
assertion  (\ref{eq.DistrH0}) of Theorem~\ref{t.asy-test} is established under a different set of weak dependence conditions on the underlying functional processes than those stated in Assumption 1.}  
\end{remark}

\begin{prop} \label{le.power} 
Suppose that the conditions of Theorem~\ref{t.asy-test} are satisfied.
Then, under $\mathcal H_1$ and  as $ T \rightarrow \infty$,
\begin{align*}
%\frac{\sqrt{b}T\, \mathcal U_T -b^{-1/2}\widehat{\mu}_0}{\widehat{\theta}_0} 
t_{\mathcal U} & =  
 \sqrt{b}T \int_{-\pi}^{\pi}\|{\mathcal F}_{X,\lambda} -{\mathcal F}_{Y,\lambda}\|^2_{HS} d\lambda + o_P(\sqrt{b} T) \ \rightarrow +\infty, \quad \mbox{in probability}.
\end{align*}
\end{prop}

The above  result, together with  Theorem~\ref{t.boo-test}  and  Slutsky's theorem,  imply that the power of the studentized test $t_{\mathcal U}$  based on the bootstrap critical values obtained from the distribution of  the bootstrap studentized test $t^\ast_{\mathcal U}$ converges to unity as $ T\rightarrow \infty$, i.e., the test  $t_{\mathcal U}$ is consistent.

\section{Numerical Results}
\subsection{Choice of the Smoothing Parameter}
Implementing the studentized test $t_{\mathcal U}$ requires the choice of the smoothing bandwidth $b$.  For
univariate and multivariate time series, this issue has been investigated in the context of a cross-validation type criterion by Beltr\~ao and Bloomfield \cite{BB-1987}, Hurvich \cite{H-1985} and Robinson \cite{R-1991}. 
However,  adaption of the multivariate approach of Robinson \cite{R-1991} to the  spectral density estimator  
$\widehat{f}_{X,\lambda}(\sigma_r,\tau_s)$, for $r,s \in \{1, \ldots, k\}$, faces problems  due to the  high dimensionality of the periodogram operator involved. 

We propose a simple approach to select the  bandwidth $b$  used in  our testing procedure which  is based on the idea to overcome the high-dimensionality of the  problem by selecting a single 
  bandwidth based on the ``on average" behavior of the pooled estimator $ 
  \widehat{f}_{\lambda}(\sigma_r,\tau_s)$, that is, its behavior  over all points  $r,s \in \{1, \ldots, k\}$ in $[0,1]^2$ for which the functional random elements $ X_t$ and $ Y_t$   are observed. 
 To elaborate, define first the following quantities. The averaged periodogram
 \[ \widehat{I}_T(\lambda) = \frac{1}{k^2}\sum_{r=1}^k\sum_{s=1}^k\big\{
  \frac{1}{2}\widehat{p}_{X,\lambda}(\sigma_r,\tau_s) + \frac{1}{2}\widehat{p}_{Y,\lambda}(\sigma_r,\tau_s)\big\}\]
  and the averaged pooled spectral density estimator 
 \[ \widehat{g}_b(\lambda) =\frac{1}{k^2}\sum_{r=1}^k\sum_{s=1}^k\big\{
  \frac{1}{2}\widehat{f}_{X,\lambda}(\sigma_r,\tau_s) + \frac{1}{2}\widehat{f}_{Y,\lambda}(\sigma_r,\tau_s)\big\}.\]
  Notice that  $\widehat{I}_T(\lambda)$ can be interpreted as  the periodogram at frequency $ \lambda$ of the pooled, real-valued  univariate process 
  $ \{V_t= \frac{1}{2}\int_0^1 X_t(s)ds + \frac{1}{2}\int_0^1 Y_t(s)ds, t \in {\mathbb Z}\}$ while $\widehat{g}_b(\lambda)$ is  an estimator of the spectral density $g$ of $ \{V_t,t\in\Z\}$. We then choose the bandwidth $b$ by minimizing the objective function 
  \[ CV(b) = \frac{1}{N} \sum_{t=1}^N \big\{ {\log(\widehat{g}_{-t}(\lambda_t))} + \widehat{I}_T(\lambda_t)/\widehat{g}_{-t}(\lambda_t)\big\},\]
over a grid of values of $b$,  where 
$\widehat{g}_{-t}(\lambda_t) = (Tb)^{-1}\sum_{s\in N_t} W((\lambda_t-\lambda_s)/b) \widehat{I}_T(\lambda_s)$ 
and $N_t=\{s: -N \leq s \leq N \ \mbox{and} \ s \neq \pm t\}$. 
That is, $ \widehat{g}_{-t}(\lambda_t)$ is the leave-one-out kernel estimator of 
$ g(\lambda)$, i.e., the estimator obtained after deleting the $t$-th frequency; see also Robinson \cite{R-1991}.\\  

Due to the computational complexity of the simulation analysis studied in the next section, the use of this automatic choice of the bandwidth $b$ will only be illustrated in the real-life data example considered in Section 5.3.

\subsection{Monte-Carlo Simulations}
{We generated functional time series stemming from the following functional moving average (FMA) processes,
\begin{align}
\label{mod-1}
X_t&=A_1(\varepsilon_{t-1}) +\alpha_2 \varepsilon_{t-2}+\varepsilon_t, \\
Y_t&=A_1(e_{t-1})  +e_t,
\label{mod-2}
\end{align}
$t \in \{1,\ldots,T\}$, where the $\varepsilon_t$  and $e_t$ are generated as independent from each other  i.i.d Brownian bridges and $A_1 $ is an integral operator
with  kernel function $\psi(\cdot,\cdot)$ given by
$$
\psi(u,v)=\frac{e^{-(u^2+v^2)/2}}{4\int_{0}^{1}e^{-t^2}dt}, \quad (u,v) \in [0,1]^2.
$$
%\begin{equation}
%X_t=\varepsilon_t + a_1 \varepsilon_{t-1}+ a_2 \varepsilon_{t-2}, 
%\label{mod-1}
%\end{equation}
%and
%\begin{equation}
%\quad Y_t=\varepsilon_t + a_1 \varepsilon_{t-1},
%\label{mod-2}
%\end{equation}
%$ t=1,2,\ldots,T$, where $\varepsilon_t$'s are i.i.d.~Brownian Bridges. 
All curves were approximated using 21 equidistant points in the unit interval and transformed into functional objects using the Fourier basis with 21 basis functions. Three sample sizes $T=50$, $T=100$ and $T=200$ were considered and the bootstrap test was applied using three nominal levels, $\alpha=0.01$,  $\alpha=0.05$ and $\alpha=0.10$. All bootstrap calculations were based on $B=1,000$ bootstrap replicates and $R=500$ model repetitions.  To investigate the empirical size and power behavior of the bootstrap test, we 
%set $a_1=0.8$ and 
consider a selection of $a_2$ values, i.e., $a_2 \in \{0.0, 0.2, 0.4, 0.6, 0.8, 1.0\}$, and various bandwidths $b$. 
%avoiding a cross validation choice of $b$ for computational time reasons. 
(Notice that $a_2=0$ corresponds to the null hypothesis while $ a_2\neq 0$ to the alternative.) 
\vspace*{0.3cm}
}

\begin{figure}[h] 
\begin{center}
\vspace*{-0.8cm}
\includegraphics[angle=0,height=6.0cm,width=7.0cm]{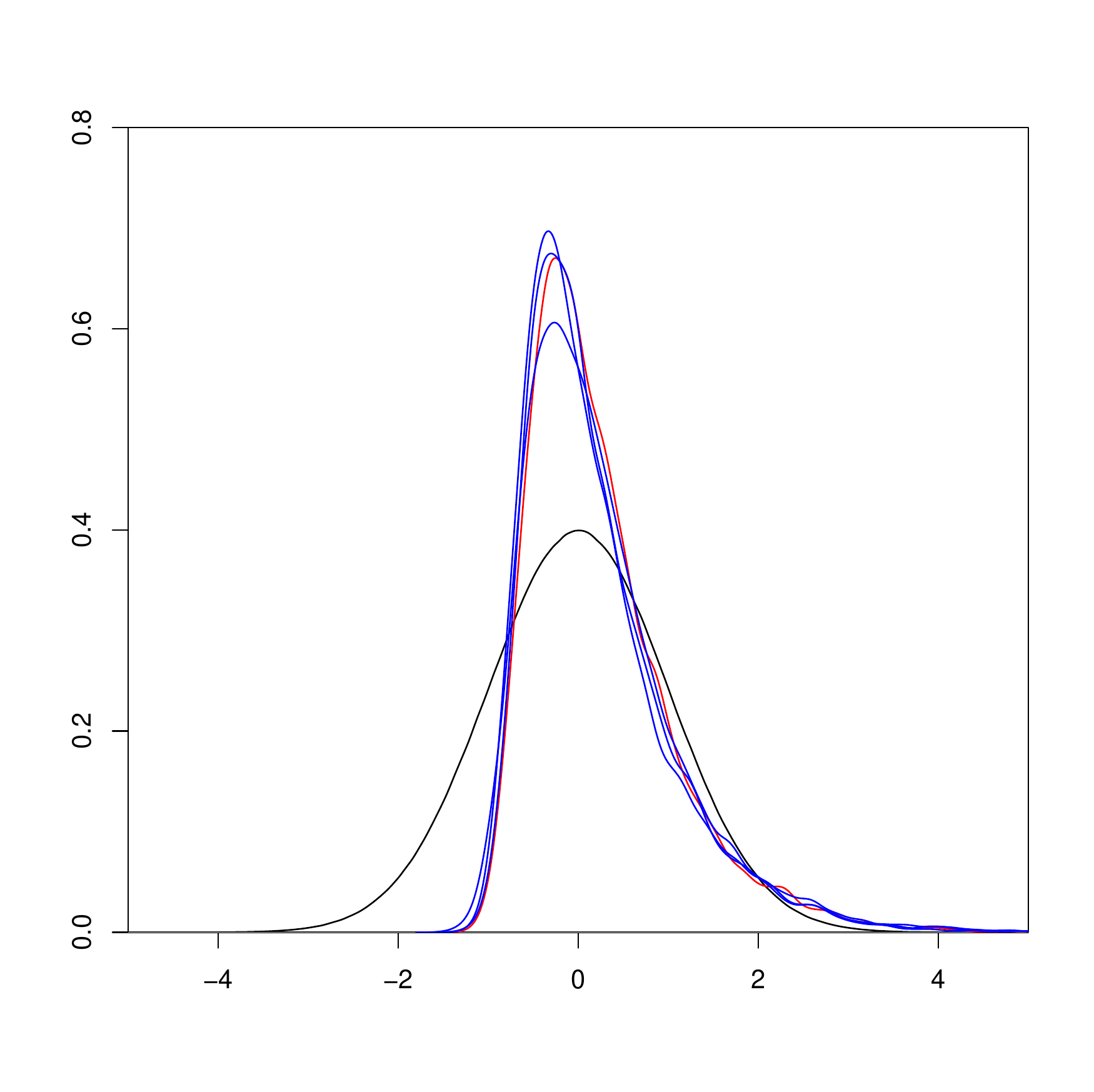}
\includegraphics[angle=0,height=6.0cm,width=7.0cm]{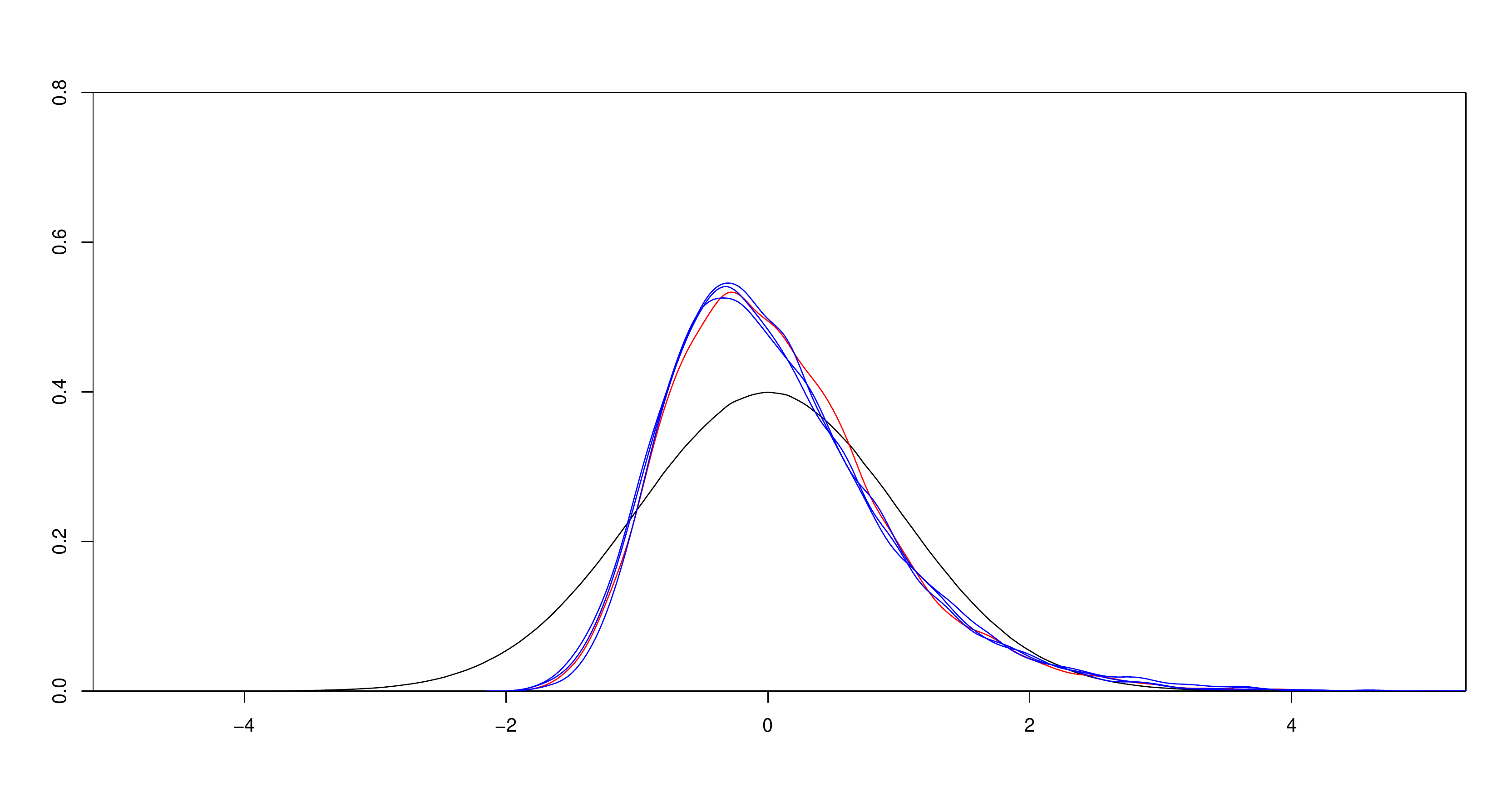}
\end{center}
\caption{Density plots of the estimated exact standardized distribution of $t_{\mathcal U}$ (red line), the standard Gaussian {distribution} (black line) and three bootstrap approximations (blue lines). Left panel, $T$=50 ($h=0.2$), right panel $T$=500 ($h=0.04$).}
\label{fig.dens}
\end{figure}

We first demonstrate the ability of the bootstrap procedure to approximate the distribution of the test statistic under the null. For this,  and in order to estimate the exact distribution of the studentized   test $t_\mathcal U$ (see (\ref{eq.stud-test})), 10,000 replications of the process  (\ref{mod-1}) and (\ref{mod-2}) with $a_2=0$ have been generated, and a kernel density estimate of this exact  distribution has been obtained using a Gaussian kernel with bandwidth $h$. The suggested bootstrap procedure is then applied to three randomly selected time series and the bootstrap studentized test $ t^\ast_\mathcal U$  (see (\ref{tst-fanis})) has been calculated. Two sample sizes of $ T=50$ and $T=500$ observations have been considered. Fig.~\ref{fig.dens} shows the results obtained together with the approximation of the distribution of $t_\mathcal U$ provided by the central limit theorem, i.e., the $\mathcal N(0,1)$ distribution. As it can be  seen from this figure, the convergence towards  the asymptotic Gaussian distribution is very slow. Even for sample sizes as large as $T=500$, the exact distribution retains its skewness which is not reproduced by the $\mathcal N(0,1)$ distribution. In contrast to this, the bootstrap approximations are very good and the estimates  of the exact densities, especially in the critical right hand tale of this distribution, are very accurate. \\

%\vspace*{0.6cm}

%\begin{center}
%Figure 1 about here. 
%\end{center}

We next investigate the finite sample size and power behavior of the bootstrap studentized test under the aforementioned  variety of  process parameters and three different sample sizes, $T=50$, $T=100$ and $T=200$.  The results obtained are shown in Table \ref{tab:empsize-power}. As it is evident from this table, the bootstrap studentized test shows a very good empirical size and power behaviour even in the case of  $T=50$ observations. In particular, the empirical sizes are close to the nominal ones and the empirical power of the test increases to one as the deviations from the null become larger (i.e., larger values of $a_2$) and/or the sample size increases.

\begin{table}[t]
\caption{Empirical size and power of the bootstrap studentized test for  functional time series generated according to models (\ref{mod-1}) and (\ref{mod-2}). \label{tab:empsize-power}}
\centering
\begin{tabular}{ccccccccc}
\hline
%\hline
 %&  &  &  & &  & \\
                   &                       &   & b=0.2& &&&b=0.3 &  \\
 %                  &              &         & &    $\alpha$     & &   \\
       $T$ & $a_2$ & $\alpha=0.01$ & $\alpha=0.05$ & $\alpha=0.10$ & & $\alpha=0.01$ & $\alpha=0.05$ & $\alpha=0.10$\\
%&  &  &  & &  & \\
\hline
%&  &  &  & &  & \\
50 &  0.0   &  0.010  &0.048  &0.096 & &0.020  &0.058&0.106 \\
         &           0.2 & 0.016 &0.082  &0.158 & &0.030  &0.092&0.164  \\ 
         &            0.4 & 0.062  &0.238  & 0.338 & &0.048 &0.154&0.276\\
         &            0.6 & 0.178  & 0.390  & 0.518 & &0.124  &0.334&0.500\\
         &            0.8 &  0.346  &0.616  & 0.736& & 0.258&0.502&0.670\\         
         &            1.0 &  0.488  &0.768  & 0.872& & 0.464&0.728&0.840\\ 
         &  &  &  & &  & & & \\ 
%\hline
%\hline        
&  &    & b=0.1&  & & & b=0.2&\\ 
      & $a_2$ & $\alpha=0.01$ & $\alpha=0.05$ & $\alpha=0.10$ & & $\alpha=0.01$ & $\alpha=0.05$ & $\alpha=0.10$\\ 
   \hline      
100 &  0.0   &  0.018  &  0.050& 0.092 & &0.008  &0.046&0.080\\
         &           0.2 & 0.028&  0.112& 0.210 &  &0.028  &0.112&0.196   \\ 
         &            0.4 &  0.138&  0.328& 0.472 & &0.122  &0.344&0.470  \\
         &            0.6 &  0.382 &   0.652& 0.764 & &0.374  &0.622& 0.766\\
         &            0.8 & 0.650&  0.858& 0.922  & &0.624  &0.836&0.922  \\          
         &            1.0 & 0.872&  0.968& 0.984  & &0.874  &0.966&0.990  \\  
%\hline
%\hline    
&  &  &  & &  & & & \\     
&  &    & b=0.06&  & & & b=0.1&\\ 
    & $a_2$ & $\alpha=0.01$ & $\alpha=0.05$ & $\alpha=0.10$ & & $\alpha=0.01$ & $\alpha=0.05$ & $\alpha=0.10$\\ 
   \hline      
200 &  0.0   &0.014  &0.042  &0.088 &&0.004  &0.044&0.100\\
         &           0.2 &0.046  &0.154  &0.272 & & 0.056 & 0.164&0.290  \\ 
         &            0.4 &0.298 &0.576 &0.698  &&  0.364&  0.620& 0.760\\
         &            0.6 &0.708  &0.910   &0.956  & &0.788&0.956&0.978\\
         &            0.8 &0.924  &0.992  &0.998& &0.960 &0.996&0.998  \\          
         &            1.0 &0.992  &1.000  &1.000& &1.000 &1.000&1.000  \\            
\hline
%\hline
\end{tabular}
%\caption{Empirical size and power of the bootstrap studentized test for  functional time series   generated according to models (\ref{mod-1}) and (\ref{mod-2}).\label{tab:empsize-power}}
\end{table}

\subsection{A Real-Life Data Example} 
\begin{figure}[t]
  \centering
 \includegraphics[trim = 10mm 15mm 10mm 15mm,clip=true,width=0.85\textwidth]{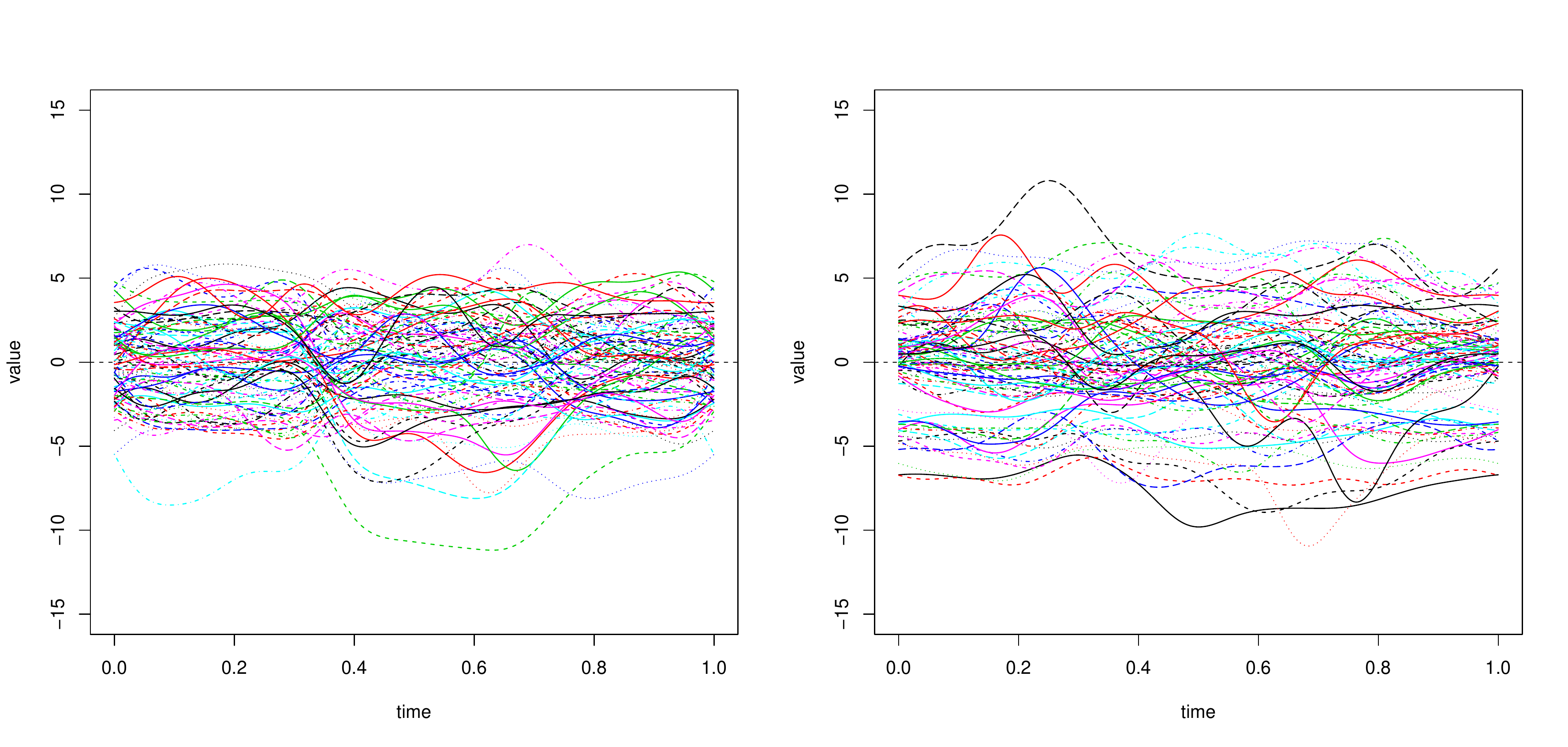}
 % \caption{1a}
  %\label{fig:sfig1}
\caption{{ Centered temperature curves} of winter period (left  panel) and of summer period (right panel). There are 92 centered curves in each period, rescaled in order to be defined in the unit interval.}
\label{fig:win-sum}
\end{figure}

We applied the bootstrap studentized test  to a data set consisting of temperature measurements recorded in Nicosia, Cyprus, for 
the winter period,  December 2006 
to beginning of March 2007 and for the summer period, June 2007 to end of August 2007. 
%In order to have the same number of days (92) in each of the two periods (i.e., winter and summer periods, respectively), the temperature measurements of the 1st and 2nd of March 2007 were included in the winter period. 
It is well-known that the mean temperatures during winter periods are smaller than those of summer periods. Our aim is to test whether there is also a significant difference in the autocovariance structure of the winter and summer periods. The data consists of two samples of curves  $\{(X_{t}, Y_t),\;t \in \{1,\ldots, 92 \}\}$, where $X_{t}$ represents the temperature of day $t$ for  Dec2006-Jan2007-Feb2007-March2007  and $Y_t$ for Jun2007-Jul2007-Aug2007. More precisely, $X_{1}$ represents the temperature of the 1st of December 2006 and $X_{92}$ the temperature of the 2nd of March 2007, whereas $Y_{1}$ represents the temperature of the 1st of June 2007 and $Y_{92}$ the temperature of the 31st of August 2007. The temperature recordings were taken in $15$ minutes intervals, i.e., there are $k=96$ temperature measurements for each day for a total of $T=92$ days in both groups. These measurements were transformed into functional objects using the Fourier basis with 21 basis functions. All curves were rescaled in order to be defined in the unit interval. Fig.~\ref{fig:win-sum} shows the  centered temperature curves of the winter and summer periods, i.e., the curves in each group are transformed  by subtracting the corresponding group sample mean functions.
\vspace*{0.3cm}

%The bootstrap $p$-values for a selection of bandwidths $b$, using $B=10,000$ bootstrap replicates, are: 0.0986 ($b=0.04$), 0.047 ($b=0.06$), 0.0255 ($b=0.08$), ?? ($b=0.10$), ?? ($b=0.12$), 0.0060 ($b=0.14$), leading to the rejection of the null hypothesis for all commonly $\alpha$-levels. 

Using the cross-validation algorithm described in Section 5.1, the bandwidth chosen is equal to  $b_{CV}=0.075$  and the 
corresponding $p$-value of the bootstrap based studentized test  is equal to 0.030 (based on $B=10,000$ bootstrap replications), leading to a rejection of the null hypothesis for almost all commonly used $\alpha$-levels.  This implies that the dependence properties, as measured by autocovariances, of the temperature measurements of the winter period differ significantly from those of the summer period.
\vspace*{0.3cm}

%\begin{figure}
%  \centering
%  \includegraphics[trim = 10mm 15mm 10mm 15mm,clip=true,width=0.65\textwidth]{cross-val}
%%  \caption{1a}
%%  \label{fig:sfig1}
%\caption{Plot of CV(b) against a grid of  values of $b$.}
%\label{fig:win-sum}
%\end{figure}

\begin{figure}
  \centering
  \includegraphics[clip=true,width=0.45 \textwidth]{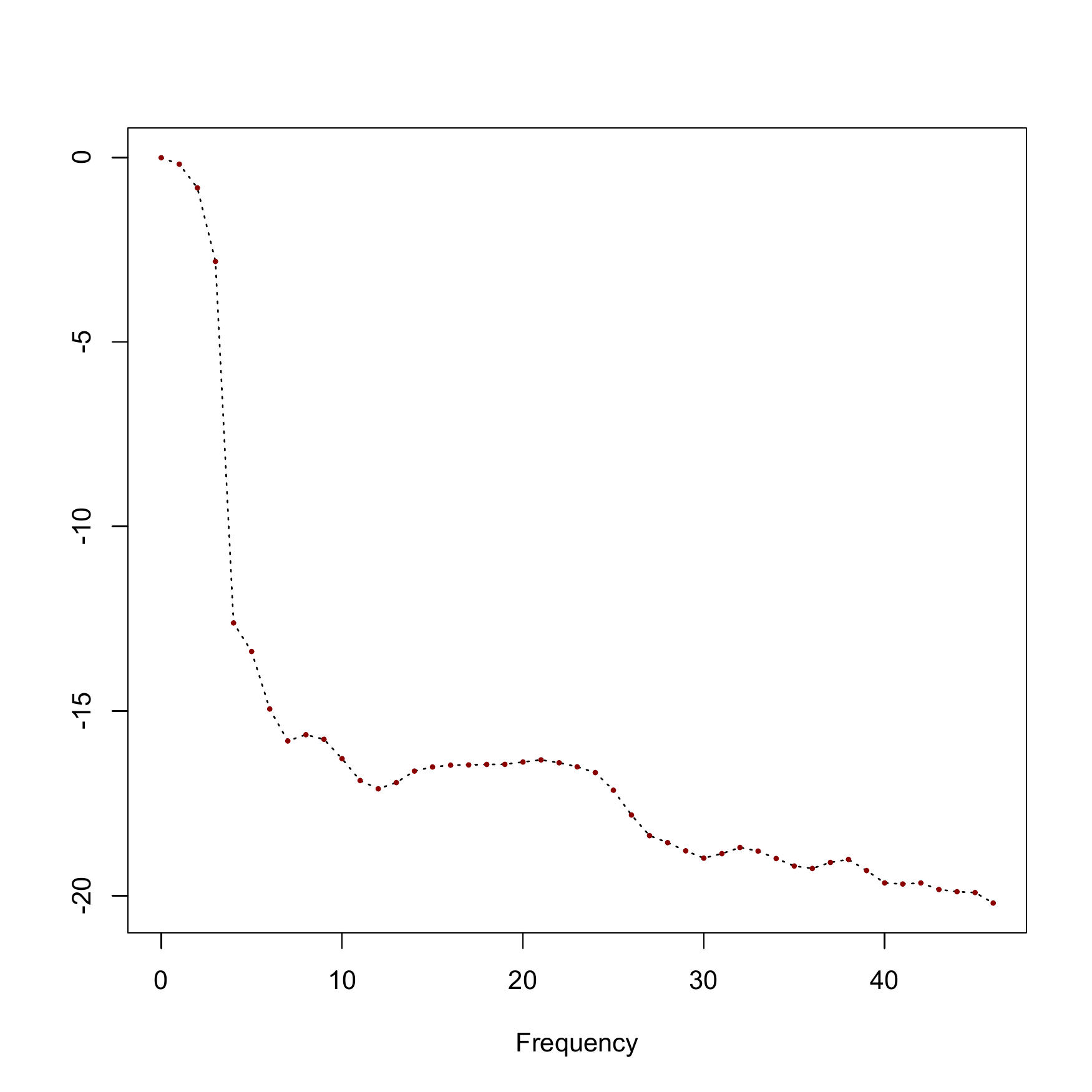}
 \includegraphics[clip=true,width=0.45 \textwidth]{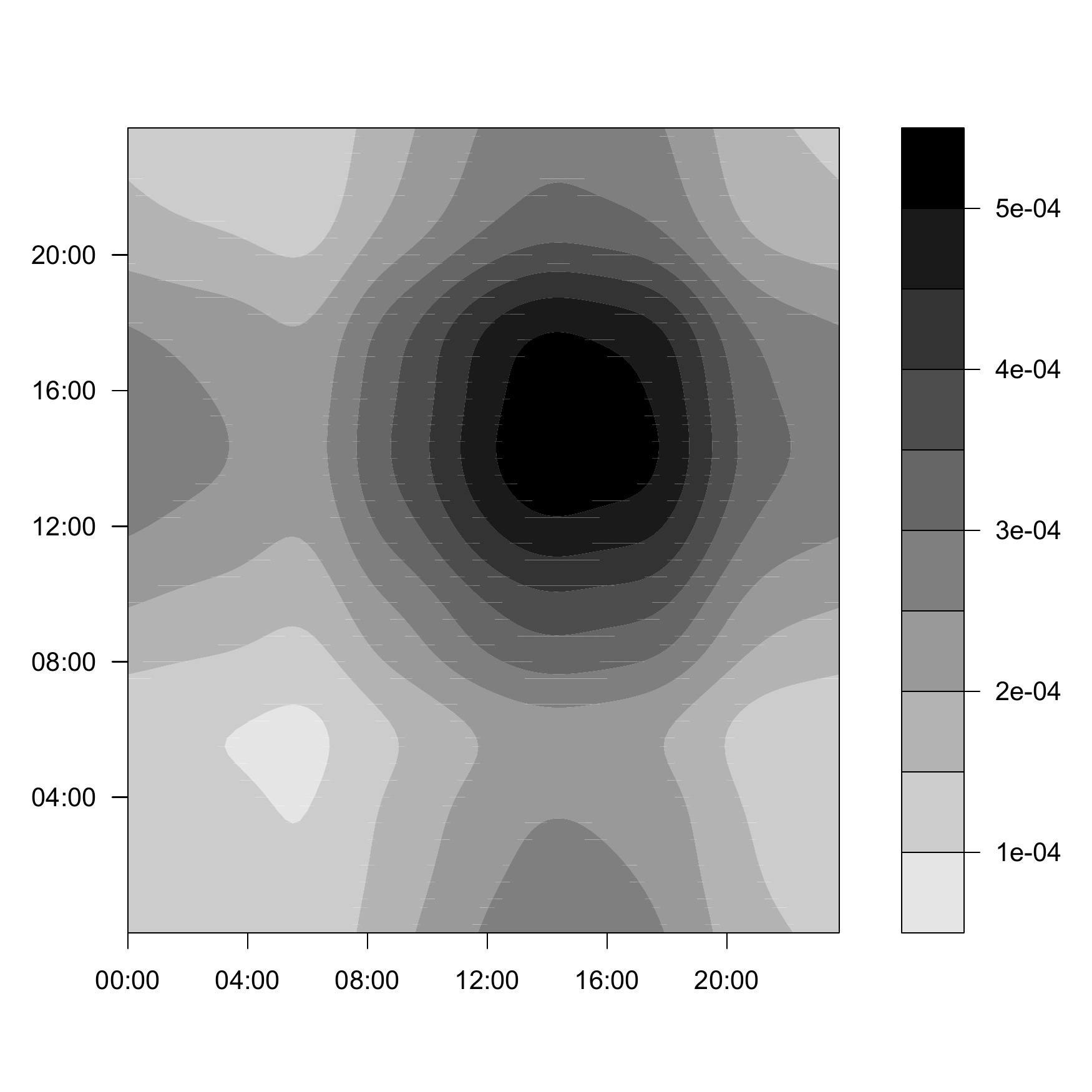}
  %\caption{1a}
  %\label{fig:sfig1}
\caption{(a)   Plot of $\widehat{Q}_{T,\lambda_j}$  (vertical axes, log-scale) against the frequencies $ \lambda_j$, ${j \in \{0,1, \ldots, N\}}$  (horizontal axes), for the temperature data, using the bandwidth $b=b_{CV}$=0.075.
(b) Plot of the difference $ \widehat{D}_T(\sigma_r,\tau_l)$ between the estimated spectral density kernels using the  same bandwidth for the values $ (\sigma_r,\tau_l)$, with $(r, l) \in \{1, \ldots, 96\}\times\{1, \ldots, 96\}$. }
\label{fig.fre-plot}
\end{figure}

In order to  understand the reasons leading to this rejection,  we decompose  the standardized test $ t_{\mathcal U}$  after  ignoring the centering sequence  $b^{-1/2} \widehat\mu_0$ and approximating  the integral of the (squared) Hilbert-Schmidt norm by the corresponding Riemann sum over the Fourier frequencies $ \lambda_j=2\pi j/T$, as follows:
\begin{align} \label{eq.dec}
\sqrt{b}T \, \mathcal U_T\Big/\widehat\theta_0 & \approx  2\pi \sqrt{b}\sum_{\substack{j=-N\\j\neq0}}^{N}  \|\widehat{\mathcal F}_{X,\lambda_j} - \widehat{\mathcal F}_{Y,\lambda_j}\|^2_{HS}\Big/\widehat\theta_0 = \sum_{\substack{j=-N\\j\neq0}}^{N} \widehat{Q}_{T,\lambda_j},
\end{align}
where
\[ \widehat{Q}_{T,\lambda_j} = 2\pi \sqrt{b} \, \|\widehat{\mathcal F}_{X,\lambda_j} - \widehat{\mathcal F}_{Y,\lambda_j}\|^2_{HS}\Big/\widehat\theta_0  \geq 0.\]
Expression (\ref{eq.dec}) shows  the contributions of the differences 
$  \|\widehat{\mathcal F}_{X,\lambda_j} - \widehat{\mathcal F}_{Y,\lambda_j}\|^2_{HS} $ for each frequency $ \lambda_j$  to the total value of the test statistic $ \mathcal U_T$.  Large values of $\widehat{Q}_{T,\lambda_j} $ pinpoint, therefore,  to frequency regions  from which   large contributions to the test statistic $ \mathcal U_T$ occur.  
%Notice that  $Q_{T,\lambda_j} $ can be estimated by $ \widehat{Q}_{T,\lambda_j} = 2\pi \sqrt{b}\,  \|\widehat{\mathcal F}_{X,\lambda_j} - 
%\widehat{\mathcal F}_{Y,\lambda_j}\|^2_{HS}\big/\widehat{\theta}_0 $, where $ \widehat{\theta}_0$ is  the  estimator of $ \theta_0$ discussed 
%in Section 3.  
 A  plot of  the estimated quantities $\widehat{Q}_{T,\lambda_j} $ against the frequencies $ \lambda_j$, $j \in \{0, \ldots, N\}$, is, therefore,  very informative   in  identifying  frequency regions where  differences between the two  spectral density operators are large and is very helpful  for  interpreting  the  results of the testing procedure. 
 
 Complementary to the decomposition $\widehat{Q}_{T,\lambda_j} $ of the test statistic $ {\mathcal U}_T$, one also can  identify the regions in $ [0,1]\times[0,1]$ 
 which deliver large contributions to the test statistic and which lead to a  rejection of the null hypothesis. In particular, the test statistic also can be written as 
 \[   \widehat{Q}_{T,\lambda_j}  \approx \sum_{r=1}^T\sum_{l=1}^T \widehat{D}^2_T(\sigma_{r},\tau_{l}) , \ \ \widehat{D}_T^2(\sigma,\tau) =\frac{2\pi \sqrt{b}}{T^2}  \sum_{\substack{j=-N\\j\neq0}}^{N}\big|\widehat{f}_{X,\lambda_j}(\sigma,\tau) - \widehat{f}_{Y,\lambda_j}(\sigma,\tau)\big|^2\Big/\widehat{\theta}_0.\]
 Notice that $ \widehat{D}_T^2(\sigma_r,\tau_l)$ shows the contribution of the differences between the estimated spectral density kernels  (averaged over all Fourier frequencies) at  the  points  $ (\sigma_r,\tau_l)\in [0,1]\times[0,1]$ to the test statistic $ {\mathcal U}_T$.  
 Large values of $ D^2_T(\sigma_{r},\tau_{l})$ pinpoint to points   $(\sigma_{r},\tau_{l} )\in [0,1]\times[0,1]$ where large differences (averaged over all frequencies) between  the corresponding spectral density kernels occur. 
 Combined with the frequency decomposition $ Q_{T,\lambda_j}$, the decomposition $ D_T(\sigma_r,\tau_l)$ may further help in better understanding the test results.

\vspace*{0.3cm}

Fig.~\ref{fig.fre-plot}(a) shows  for the real-life temperature data example considered the  plot of $ \widehat{Q}_{T, \lambda_j}$ at a log-scale.  Fig.~\ref{fig.fre-plot}(b) shows,  for the same data set,  a plot of 
the differences $ \widehat{D}_T(\sigma_r,\tau_l)$. 
%between the estimated spectral density   kernels $ |\widehat{f}_{X,\lambda_j}(\sigma_r,\tau_l) - \widehat{f}_{Y,\lambda_j}(\sigma_r,\tau_l) |^2 $ averaged  over all Fourier frequencies $ \lambda_j$.  
As it can be
seen from Fig.~\ref{fig.fre-plot}(a), the large values of the test statistic $ \mathcal U_T$ which leads to a rejection of the null hypothesis, are mainly due to the  large differences between the two spectral density operators at  the low frequency region. That  is, differences in the long term periodicities  between  the winter and  the summer temperature curves seem to be the main reason for rejecting the null hypothesis. Fig.~\ref{fig.fre-plot}(b) shows that the main differences  between the spectral density kernels of the two functional time series, occur in the  afternoon period and, more specifically, between the hours  12.00 to  4.00 p.m.  The differencies of the (averaged) spectral density kernels,  for values of $ \tau$ and $\sigma$ within this time frame, seem to be the largest.  
These findings are  probably due to the fact that in Cyprus, compared to the rather day-long stable weather conditions of the summer period, the weather conditions in the winter period are more volatile, change gradually during the day and reach their peak in the afternoon. 

%
%\begin{center}
%Figure 2 about here. 
%\end{center}
%%
%\begin{center}
%Figure 3 about here. 
%\end{center}

%\section{Conclusions} {We have proposed  a fully functional,  frequency domain $L^2$-type test for  testing  equality of the entire second order structure of two independent  functional time series. A large sample  Gaussian approximation of the distribution of  the  test   statistic under the null has been derived.  Moreover, in order to improve upon the  Gaussian approximation, a 
%frequency  domain bootstrap procedure has  been proposed which leads to more  accurate estimates of the distribution of the test statistic of interest under the null. 
%Procedures to implement the test statistic under the null and to select the bandwidth parameters involved have also been presented. 
%The  finite sample behavior of the test has been investigated by means of simulations and a real-life  data example has been analyzed using the method presented in this paper. }

%\appendix

%%%%%%%%%%%%%%%%%%%%%%%%%%%%%%%%%%%%%%%%%%%%%%%%%%%%%%%%%%%%%%%%%%%%
%\section{Auxiliary Results and Proofs}
%%%%%%%%%%%%%%%%%%%%%%%%%%%%%%%%%%%%%%%%%%%%%%%%%%%%%%%%%%%%%%%%%%%%

\section*{Acknowledgments}
We thank the Editor-in-Chief, the Managing Editor of the SI JMVA FDA  and the two referees for their helpful comments. This research is partly funded by the Volkswagen Foundation (Professorinnen f\"ur
Niedersachsen des Nieders\"achsischen Vorab) and by a University of Cyprus Research Grant.

\section*{Appendix: Auxiliary Results and Proofs}

First, we introduce some notation that will be used throughout our proofs. 
$ \|\cdot\|_2$ denotes the norm of $ L^2$,  $ \|\cdot\|_{\mathcal N}$ the nuclear norm of an operator $T$,  $T^*$ is the adjoint operator and $ \langle\cdot, \cdot \rangle_{HS}$ the inner product on the space of Hilbert-Schmidt operators; see the Supplementary Material for more details.  Furthermore, we write $A(e^{-i\lambda})=\sum_{j\in\Z}A_je^{-ij\lambda}$ with  the operators $A_j$ defined as in Assumption~1. The periodogram operators of the innovations time series
$ \varepsilon_t$ and $ e_t$, $t=1, \ldots, n$, at frequency $\lambda$,  $I_{\varepsilon,\lambda}$ and $I_{e,\lambda}$, respectively, are defined as the integral operators induced by right integration of 
\begin{equation}\label{eq.p.center}
\begin{aligned}
\widehat p_{\varepsilon, \lambda}(\sigma,\tau)&=\frac{1}{2\pi T}\sum_{s_1,s_2=1}^T \varepsilon_{s_1}(\sigma)\varepsilon_{s_2}(\tau)\, \exp(-i\lambda( s_1- s_2)),\\ 
\widehat p_{e, \lambda}(\sigma,\tau)&=\frac{1}{2\pi T}\sum_{s_1,s_2=1}^T e_{s_1}(\sigma)e_{s_2}(\tau)\, \exp(-i\lambda( s_1- s_2)).
\end{aligned}
\end{equation}
The centered counterparts are denoted by $I_{\varepsilon,\lambda}^c$ and $I_{e,\lambda}^c$. Finally, define
$$
Q_{X,\lambda}^{c} :=A(e^{-i\lambda}) I^c_{\varepsilon,\lambda} A(e^{-i\lambda})^*, \ \ 
Q_{Y,\lambda}^{c} :=B(e^{-i\lambda}) I^c_{e,\lambda} B(e^{-i\lambda})^* 
.
$$
Here, $S T$ denotes the composition $S(T(\cdot))$ of the operators $S$ and $T$.
\vspace*{0.3cm}

\noindent
{\bf Proof of Lemma \ref{l.ex-SDO}.}
The assertions of the lemma are immediate consequences of Proposition~2.1 in Panaretos and Tavakoli \cite{PT-2013} if 
$
\sum_{t\in\Z} \|\mathcal R_{X,t}\|_{\mathcal N}<\infty $ and $
\sum_{t\in\Z} \|r_{X,t}\|_{2}<\infty
$
 and similar results for the process $(Y_t)_{ t \in \Z}$ hold true.
The first inequality follows from  expression (1.4)  of the supplement.
For the second result, use the expression ${\mathcal R}_{X,t}=\sum_{j\in \Z} A_{j+t} C_\varepsilon A^\ast_j $, see the Supplement Material,  and get  
%we apply \textbf{Lemma 3.3 of Rademacher}.
%$$
%\begin{aligned}
$\sum_{t\in\Z}\|r_{X,t}\|_{2}=\sum_{t\in\Z}\|\mathcal R_{X,t}\|_{HS}
%&=\sum_{t\in\Z} \left\|\sum_{j\in\Z} A_{j+t} \mathcal R_{\varepsilon,0} A_j^*\right\|_{HS}\\
\leq \sum_{t\in\Z}\sum_{j\in\Z} \|A_{j+t}\|_{\mathcal L}\,\| C_{\varepsilon}
\|_{HS}\, \|A_{j}\|_{\mathcal L}$,
%\end{aligned}
%$$
which is finite under Assumption~1. 
\vspace*{0.3cm}

The proof of Theorem~\ref{t.asy-test}  uses the following two  lemmas, the proofs of which are given in the Supplementary Material.

\begin{lem}\label{l.mean}
Suppose that the assumptions of Theorem~\ref{t.asy-test} hold true. Then
$$
\sqrt bT\, M_{T,0}-b^{-1/2} \, \mu_0=o_P(1),
$$
where
\begin{equation}\label{eq.MT}
\begin{aligned}
M_{T,0}&=\int_{-\pi}^\pi\frac{1}{b^2T^2}\sum_{t=-N}^N W^2\left(\frac{\lambda-\lambda_t}{b}\right)\left\|Q_{X,\lambda_t}^{c}-Q_{Y,\lambda_t}^{c}\right\|_{HS}^2d\lambda \,.
\end{aligned}
\end{equation}
\end{lem}
\vspace*{0.3cm}

\begin{lem}\label{l.var}
Suppose that the assumptions of Theorem~\ref{t.asy-test} hold true. Then, 
$$
{\rm var}(\sqrt b T\, L_{T,0})\ninfty\theta_0^2
$$
for $\theta_0$ defined in Theorem~\ref{t.asy-test}, where 
\begin{equation}\label{eq.LT}
\begin{aligned}
L_{T,0}&= \frac{1}{b^{2}T^2}\int_{-\pi}^\pi\int_0^1\int_0^1\sum_{t_1,t_2=-N,\, t_1\neq t_2}^N \, W\left(\frac{\lambda-\lambda_{t_1}}{b}\right)\,W\left(\frac{\lambda-\lambda_{t_2}}{b}\right)\, d\lambda
\\
&\quad\times\Big\langle Q_{X,\lambda_{t_1}}^{c}-Q_{Y,\lambda_{t_1}}^{c}\;, Q_{X,\lambda_{t_2}}^{c}- Q_{Y,\lambda_{t_2}}^{c}\Big\rangle_{HS}.
\end{aligned}
\end{equation}
\end{lem}
\vspace*{0.3cm}

\noindent
{\bf Proof of Theorem~\ref{t.asy-test}.}
From Theorem 1.2 of the Supplementary  Material we obtain
$$
I_{X,\lambda}= A(e^{-i\lambda}) I_{\varepsilon,\lambda} A(e^{-i\lambda})^*+R_{T,\lambda}\quad\text{with}\quad \sup_{\lambda\in\{2\pi t/T|t=-N,\dots,N\}} {\rm E}\|R_{T,\lambda}\|_{HS}^2=O(T^{-1}).
$$
 This gives
\begin{align}
\label{eq.clt1}
\sqrt b T \,\mathcal U_T &=\sqrt b T \int_{-\pi}^\pi\left\|\frac{1}{bT}\sum_{t=-N}^N W\left(\frac{\lambda-\lambda_t}{b}\right)\left[A(e^{-i\lambda_t}) I_{\varepsilon,\lambda_t} A(e^{-i\lambda_t})^*-B(e^{-i\lambda_t}) I_{e,\lambda_t} B(e^{-i\lambda_t})^*\right]\right\|_{HS}^2d\lambda \nonumber \\
&\quad\;+O_P(b^{1/4}) =:\sqrt b T \,U_{T,0}+o_P(1)
\end{align}
if we can show that $\sqrt b T\, U_{T,0}=O_P(1)$. To this end, first note that under $\mathcal H_0$  $\mathcal F_X=\mathcal F_Y$. 
Now, it follows from (1.5) in the Supplementary Material   that 
$$
\frac{1}{2\pi}\,A(e^{-i\lambda_t})\, C_{\varepsilon}\,A(e^{-i\lambda_t})^*=\mathcal F_X=\mathcal F_Y=\frac{1}{2\pi}\,B(e^{-i\lambda_t})\, C_{e} \, B(e^{-i\lambda_t})^*.
$$
Additionally, we have  ${\rm E}\widehat p_{\varepsilon, \lambda}(\sigma,\tau)=c_{\varepsilon}(\sigma,\tau)/(2\pi)$  and ${\rm E}\widehat p_{e, \lambda}(\sigma,\tau)=c_{e}(\sigma,\tau)/(2\pi)$  in $L^2$ for the i.i.d.~noises. Combining both facts,  we can rewrite $U_{T,0}$ as
\begin{align}
\label{eq.clt2}
U_{T,0} =
\int_{-\pi}^\pi\left\|\frac{1}{bT}\sum_{t=-N}^N W\left(\frac{\lambda-\lambda_t}{b}\right)\left[A(e^{-i\lambda_t}) I^c_{\varepsilon,\lambda_t} A(e^{-i\lambda_t})^*-B(e^{-i\lambda_t}) I^c_{e,\lambda_t} B(e^{-i\lambda_t})^*\right]\right\|_{HS}^2d\lambda.
\end{align}
We can further split up
$
 U_{T,0}=M_{T,0}+L_{T,0}
$
where $M_{T,0}$ and $L_{T,0}$ are defined as in Lemma~\ref{l.mean} and Lemma~\ref{l.var}, respectively. In view of Lemma~\ref{l.mean}, it remains to show  that $\sqrt b T L_{T,0}\stackrel{d}{\longrightarrow} Z.$
To this end, we abbreviate
$$
w_{t_1, t_2, T}=\frac{1}{b^{3/2}T}\int_{-\pi}^\pi  W\left(\frac{\lambda-\lambda_{t_1}}{b}\right)\, W\left(\frac{\lambda-\lambda_{t_2}}{b}\right)\, d\lambda
$$
and use the Karhunen-L{\'o}eve expansion for the Gaussian innovations $\varepsilon_t$ and $e_t$. In particular, 
we have
% i.e.~there is an ONS of continuous functions $(\varphi_l)_l$ in $L_2[0,1]$  denote such that
$$
\varepsilon_s(\sigma)=\sum_l \xi_l^{(s)} \varphi_l(\sigma), \quad s\in \Z,~\sigma\in[0,1],
$$
 where $  \varphi_l \in L^2$, $ l\in \N$,   denotes  the set of  orthonormal eigenfunctions  of the operator  $C_\varepsilon$ and    the random variables $\xi_l^{(s)}=\int_0^1 \varepsilon_s(\sigma)\, \varphi_l(\sigma)\,d\sigma$ are centered normal and satisfy ${\rm cov}(\xi_{l_1}^{(s)}, \xi_{l_2}^{(s)})=0$ for $ l_1\neq l_2.$  Notice that 
 the above expression for $\varepsilon_s(\sigma)$  is valid in  $L^2$-sense and that Fubini's theorem gives ${\rm cov}(\xi_{l_1}^{(s_1)}, \xi_{l_2}^{(s_2)})=0$ for $s_1\neq s_2$.
A similar expansion holds true for $e_s$ with a possibly different set of othonormal eigenfunctions $(\phi_l)_{l\in\N}$ instead of $(\varphi_l)_{l\in\N}$. Now, we define approximating periodogram operators $I_{\varepsilon,\lambda_t}^{c,K},\; K\in \N$, with kernels 
$$
\widehat p_{\varepsilon,\lambda_t}^{ c,K}=\sum_{l_1,l_2=1}^K\varphi_{l_1}\varphi_{l_2}\frac{1}{2\pi T}\sum_{s_1,s_2}e^{i\lambda_t(s_1-s_2)}\, [\xi_{l_1}^{(s_1)}\xi_{l_2}^{(s_2)}-E(\xi_{l_1}^{(s_1)}\xi_{l_2}^{(s_2)})]
$$
and similarly for $I_{e,\lambda_t}^{c,K}$. Moreover, define
$$
Q_{X,\lambda}^{c,K} :=A(e^{-i\lambda}) I^{c,K}_{\varepsilon,\lambda} A(e^{-i\lambda})^* 
\quad\text{and}\quad  
Q_{Y,\lambda}^{c,K} :=B(e^{-i\lambda}) I^{c,K}_{e,\lambda}B(e^{-i\lambda})^* .
$$
Then, we can introduce
$$
\begin{aligned}
\sqrt b T\, L_{T,0}^{(K)}&= \sum_{t_1,t_2=-N,\, t_1\neq t_2}^N \, w_{t_1, t_2, T}
\;\Big\langle Q_{X,\lambda_{t_1}}^{c,K}-Q_{Y,\lambda_{t_1}}^{c,K}\;,Q_{X,\lambda_{t_2}}^{c,K}-Q_{Y,\lambda_{t_2}}^{c,K}\Big\rangle_{HS} =:\sum_{t_1,t_2=-N, t_1\neq t_2}^N H_{t_1,t_2,T}.
\end{aligned}
$$
From this, we get  
\begin{equation}\label{eq:var-cut}
\lim_{K\to\infty}\limsup_{T\to\infty} {\rm E}(\sqrt bT (\, L_{T,0}-\, L_{T,0}^{(K)}))^2=0.
\end{equation}
To this end, first note that under Gaussianity $|{\rm E}L_{T,0} |+|{\rm E}L_{T,0}^{(K)}|=o(1)$ for any $K$  due to independence of the spectral density operators at different frequencies $|t_1|\neq |t_2|$. Thus, it suffices to investigate ${\rm var}(\sqrt bT (\, L_{T,0}-\, L_{T,0}^{(K)}))$. With the same arguments as in the proof of Lemma~\ref{l.var} it suffices to show that
$$
\begin{aligned}
\sup_{|t_1|\neq |t_2|,|s_1|\neq |s_2|}{\rm cov}&\Big(\Big\langle Q_{X,\lambda_{t_1}}^{c}-Q_{Y,\lambda_{t_1}}^{c}\;,Q_{X,\lambda_{t_2}}^{c}-Q_{Y,\lambda_{t_2}}^{c}\Big\rangle_{HS}
  - \Big\langle Q_{X,\lambda_{t_1}}^{c,K}-Q_{Y,\lambda_{t_1}}^{c,K}\;,Q_{X,\lambda_{t_2}}^{c,K}-Q_{Y,\lambda_{t_2}}^{c,K}\Big\rangle_{HS},\\
	&\Big\langle Q_{X,\lambda_{s_1}}^{c}-Q_{Y,\lambda_{s_1}}^{c}\;,Q_{X,\lambda_{s_2}}^{c}-
	Q_{Y,\lambda_{s_2}}^{c}\Big\rangle_{HS} 
	-\Big\langle Q_{X,\lambda_{s_1}}^{c,K}-Q_{Y,\lambda_{s_1}}^{c,K}\;,Q_{X,\lambda_{s_2}}^{c,K}-Q_{Y,\lambda_{s_2}}^{c,K}\Big\rangle_{HS}\Big)
\end{aligned}
$$
converges to zero as $K\to\infty$ in the cases $t_1=\pm s_1, t_2=\pm s_2$ and $t_1=\pm s_2, t_2=\pm s_1$. Exemplarily we only  investigate
$$
\begin{aligned}
\sup_{|t_1|\neq |t_2|}{\rm cov}&\Big(\Big\langle Q_{X,\lambda_{t_1}}^{c}\;,Q_{X,\lambda_{t_2}}^{c}\Big\rangle_{HS}-\Big\langle Q_{X,\lambda_{t_1}}^{c,K}\;,Q_{X,\lambda_{t_2}}^{c,K}\Big\rangle_{HS},\\
&\Big\langle Q_{X,\lambda_{t_1}}^{c}\;,Q_{X,\lambda_{t_2}}^{c}\Big\rangle_{HS}-\Big\langle Q_{X,\lambda_{t_1}}^{c,K}\;,Q_{X,\lambda_{t_2}}^{c,K}\Big\rangle_{HS}\Big)
\end{aligned}
$$
in detail. With similar arguments as in Lemma~\ref{l.var} it can be shown that all remaining summands vanish, too. Using symmetry arguments and adding zeros, it suffices to consider\
%footnote{\textcolor{red}{DR: How does the mixed innerproduct $\langle Q_{X,\lambda_{t_1}}^{c,K}, Q_{X,\lambda_{t_2}}^c  \rangle$ come into play?- AL: OK so?}}
\begin{equation}\label{eq.sup}
\begin{aligned}
\sup_{|t_1|\neq |t_2|}{\rm cov}&\Big(\Big\langle Q_{X,\lambda_{t_1}}^{c}-Q_{X,\lambda_{t_1}}^{c,K}\;,Q_{X,\lambda_{t_2}}^{c}\Big\rangle_{HS},\Big\langle Q_{X,\lambda_{t_1}}^{c}\;,Q_{X,\lambda_{t_2}}^{c}\Big\rangle_{HS}\Big)
\end{aligned}
\end{equation}
%\textcolor{blue}{
and similar terms. To this end, let 
$$C_\epsilon^{(K)}={\rm E} \left[\left(\sum_{l=1}^K \xi_l^{(0)} \varphi_l\right)\otimes \left(\sum_{l=1}^K \xi_l^{(0)} \varphi_l\right)\right].
$$
In analogy to the proof of Lemma~\ref{l.var}, \eqref{eq.sup} can be bounded from above by
$$
\begin{aligned}
&\sup_{|t_1|\neq |t_2|} \|A(e^{-i\lambda_{t_1}})\|_{\mathcal L}^4\,\left\|{\rm E}\left([I_{\varepsilon,\lambda_{t_1}}^c - I_{\varepsilon,\lambda_{t_1}}^{c,K} ]\otimes I_{\varepsilon,\lambda_{t_1}}^c\right) \right\|_{HS}\,
\|\mathcal F_{X,\lambda_{t_2}}\|_{HS}^2 \leq \mathcal K\,\|C_\varepsilon-C_\epsilon^{(K)}\|_{HS} +o(1)
\end{aligned}
$$
for some finite constant $\mathcal K$, where the last inequality can be obtained similarly to Lemma 1.7 and Theorem 1.3 in the supplement. Mercer's Theorem finally gives $\|C_\varepsilon-C_\epsilon^{(K)}\|_{HS}\to 0$ as $K\to\infty$.
%}
We aim at applying a CLT of de Jong \cite{dJ-1987} for weighted $U$-statistics of independent random vectors. To this end,  we  rewrite 
$$
\sqrt b T\, L_{T,0}^{(K)}=\sum_{t_1,t_2=1,t_1\neq t_2}^N \widetilde H_{t_1,t_2,T}+ \sum_{t=-N}^N [H_{t,0,T}+H_{0,t,T}]+\sum_{t_1=-N}^N  H_{t_1,-t_1,T}-2H_{0,0,T},
$$
where
$$
\widetilde H_{t_1,t_2}= H_{t_1,t_2,T}+H_{-t_1,t_2,T}+H_{t_1,-t_2,T}+H_{-t_1,-t_2,T}.
$$
Straightforward calculations yield that 
$$
\sum_{t=-N}^N [H_{t,0,T}+H_{0,t,T}]+\sum_{t_1=-N}^N  H_{t_1,-t_1,T}-2H_{0,0,T}=o_P(1)
$$ in $L^2$.
Now, we apply Theorem~2.1 of de Jong \cite{dJ-1987} to
$$
\widetilde W_T=\sum_{\substack{t_1, t_2=1 \\ t_1\neq t_2}}^N \, \widetilde H_{t_1,t_2}=\sum_{\substack{t_1, t_2=1 \\ t_1\neq t_2}}^N \, \widetilde H_{t_1,t_2}(\mathbb X_{t_1},\mathbb X_{t_2}),
$$ 
where $H_{t_1,t_2}$ is a Borel function and
$$
\mathbb X_t=\frac{1}{\sqrt{2\pi T}}\sum_{s=1}^T(\xi_s^{(1)}\cos(\lambda_t s), \xi_s^{(1)} \sin (\lambda_t s),\ldots, \xi_s^{(K)}\cos(\lambda_t s), \xi_s^{(K)} \sin (\lambda_t s))'
$$ 
in their notation. 
First, note that  the assumption of Gaussian  innovations implies independence of $\mathbb X_{1},\dots, \mathbb X_N$. Moreover, this yields ${\rm E}(\widetilde H_{t_1,t_2}\mid \mathbb X_{t_1})={\rm E}(\widetilde H_{t_1,t_2}\mid \mathbb X_{t_2})=0~a.s.$ for $t_1\neq t_2$  which implies that $\widetilde W_T$ is clean (see Definition 2.1 in de Jong \cite{dJ-1987}). It remains to check conditions (a) and (b) of Theorem~2.1 of de Jong \cite{dJ-1987}. Similar to Lemma~\ref{l.var} we obtain that ${\rm var}(\widetilde W_T)$ converges to the finite constant
$$\theta_K:=\frac{4}{\pi^2}\int_{-2\pi}^{2\pi}\left\{\int_{-\pi}^\pi W(u)W(u-x)\, du\right\}^2dx\, 
\int_{-\pi}^\pi \| A(e^{-i\lambda_{t_1}})\, E[I_{\varepsilon,\lambda_{t_1}}^{c,K}]\, A(e^{-i\lambda_{t_1}})^*  \|_{HS}^4 \,d\lambda.
$$
 Subsequently, we only consider the non-trivial case of $\theta_L>0$. For condition (a), it remains to verify that
$$
\max_{t_1\in\{1,\dots, N\}}\sum_{\substack{t_2=1 \\ t_2\neq t_1}}^N {\rm var}\left(\widetilde H_{t_1,t_2}\right)=o(1).
$$  
This is an immediate consequence of ${\rm var}(H_{t_1,t_2})=0$ for $|t_1-t_2|>bT$ and
$$
\begin{aligned}
&{\rm var}(H_{t_1,t_2})= O\left(\frac{1}{b\,T^2}\right)\,=o\left(\frac{1}{bT}\right)
\end{aligned}
$$
for $|t_1-t_2|\leq bT$. 
Finally, we have to check assumption (b) of Theorem~2.1 of de Jong \cite{dJ-1987}, i.e.,
$$
{\rm E}\widetilde W_T^4\ninfty 3\theta_K^2.
$$
To this end, we argue that ${\rm E}\widetilde W_T^2\ninfty\theta_K^2$ and that the forth-order cumulant of $\widetilde W_T$ vanishes asymptotically due to the independence of the periodograms at different Fourier frequencies.
Finally, note that $\theta_K\to \theta_0$ as $K\to\infty$ which finishes the proof by Proposition 6.3.9 in Brockwell and Davis \cite{BD-1991}.\\
%\end{proof}

\noindent
{\bf Proof of Theorem~\ref{t.boo-test}.}
Recall first that in the following  calculations  all indices in the sums considered, run in the set $ \{-N, -N+1, \ldots, -1,1,\ldots, N-1, N\}$, where $ N=[(T-1)/2]$.
%; i.e., we ignore the asymptotically negligible zero frequency  in the consideration of $ \widehat{f}^\ast_{X,\lambda_t}$ and $\widehat{f}^\ast_{Y,\lambda_t} $ in Step 3 of the bootstrap algorithm. 
 Let $ \{v_j, j\in \N\} $ be an orthonormal basis  of $ L_\C^2:=L_\C^2([0,1],\mu)$
 and recall that $ \{v_i\otimes v_j, i,j\in \N\}$ is  an orthonormal basis of the Hilbert space   $HS( L_\C^2)$. The bootstrap test statistic 
\begin{equation}
\mathcal U^\ast_T = \frac{2\pi }{T}\sum_{l=-N}^N  \| \widehat{\mathcal F}_{X,\lambda_t}^\ast -  \widehat{\mathcal F}_{Y,\lambda_t}^\ast\|^2_{HS}
\end{equation}
can then be decomposed as  
\begin{align*}
 \mathcal U^\ast_T & = \frac{2\pi}{T^3b^{2}}\sum_{t_1=-N}^N\sum_{t_2=-N}^N \sum_{l=-N}^N W\Big(\frac{\lambda_l-\lambda_{t_1}}{b} \Big)
W\Big(\frac{\lambda_l-\lambda_{t_2}}{b} \Big) \langle I^\ast_{X,\lambda_{t_1}} - I^\ast_{Y,\lambda_{t_1}},  I^\ast_{X,\lambda_{t_2}} - I^\ast_{Y,\lambda_{t_2}}\rangle _{HS}\\
&= \frac{2\pi}{T^3b^{2}}\sum_{t=-N}^N\sum_{l=-N}^N W^2\Big(\frac{\lambda_l-\lambda_{t}}{b} \Big)
 \| I^\ast_{X,\lambda_{t}} - I^\ast_{Y,\lambda_{t}} \|_{HS}^2\\
& \ \ \ \ + 
\frac{2\pi}{T^3b^{2}}\sum_{\substack{t_1,t_2=-N \\ t_|\neq t_2}}^N W\Big(\frac{\lambda_l-\lambda_{t_1}}{b} \Big)
W\Big(\frac{\lambda_l-\lambda_{t_2}}{b} \Big) \langle I^\ast_{X,\lambda_{t_1}} - I^\ast_{Y,\lambda_{t_1}},  I^\ast_{X,\lambda_{t_2}} - I^\ast_{Y,\lambda_{t_2}}\rangle _{HS}\ := M^\ast_T + L^\ast_T,  
\end{align*}
with an obvious notation for  $M^\ast_T $ and $ L^\ast_T$.
In the following we use the notation 
\[ D^*_{t}(j_1, j_2) := \langle I^\ast_{X,\lambda_t} - I^\ast_{Y,\lambda_t} ,v_{j_1}\otimes v_{j_2} \rangle =   \langle J^\ast_{X,\lambda_t},v_{j_1}\rangle \langle v_{j_2}, \overline{J}^\ast_{X,\lambda_t}\rangle 
- \langle J^\ast_{Y,\lambda_t},v_{j_1}\rangle \langle v_{j_2}, \overline{J}^\ast_{Y,\lambda_t}\rangle ,\]
and the expansion 
\begin{align*}
 I^\ast_{X,\lambda_t}-I^\ast_{Y,\lambda_t} & =  J^\ast_{X,\lambda_t}\otimes \overline{J}^\ast_{X,\lambda_t} - J^\ast_{Y,\lambda_t}\otimes \overline{J}^\ast_{Y,\lambda_t} = \sum_{j_1=1}^\infty\sum_{j_2=1}^\infty  D_{t}^{{*}}(j_1,j_2) (v_{j_1}\otimes v_{j_2}).
\end{align*}
Notice that 
%using the bootstrap pseudo periodogram operator, 
$\langle J^\ast_{X,\lambda_t},v_{j}\rangle $ is for every $ j \in \N$, a complex  Gaussian random variable.
%\begin{align*}
%\mathcal U^\ast_T & = \frac{2\pi}{T^3b^2}\sum_{t_1=-N}^N\sum_{t_2=-N}^N \sum_{l=-N}^N W\Big(\frac{\lambda_l-\lambda_{t_1}}{b} \Big)
%W\Big(\frac{\lambda_l-\lambda_{t_2}}{b} \Big) \int_0^1\int_0^1 d^*_{t_1}(\sigma,\tau)\overline{d^*_{t_2}}(\sigma,\tau)d\sigma d\tau \\
%& =  \frac{2\pi}{T^3b^2}\sum_{l=-N}^N \sum_{t=-N}^N W^2\Big(\frac{\lambda_l-\lambda_{t}}{b} \Big)
% \int_0^1\int_0^1 |d^{*}_{t}(\sigma,\tau)|^2d\sigma d\tau\\
% & \ \ \ \ +  \frac{2\pi}{T^3b^2}\sum_{\substack{t_1,t_2=-N \\t_1\neq t_2}}^N \sum_{l=-N}^N W\Big(\frac{\lambda_l-\lambda_{t_1}}{b} \Big)
%W\Big(\frac{\lambda_l-\lambda_{t_2}}{b} \Big) \int_0^1\int_0^1 d^*_{t_1}(\sigma,\tau)\overline{d^*_{t_2}}(\sigma,\tau)d\sigma d\tau\\
%& = M^\ast_T + L^\ast_T,  
%\end{align*}
%with an obvious notation for $ M_T^\ast$ and $ L^\ast_T$.
We show that 
\begin{equation} \label{eq.M}
\sqrt{b} T M_T^\ast \,-\, b^{-1/2}\,\widetilde{\mu}_0 \stackrel{P}{\rightarrow} 0,
\end{equation}
and 
\begin{equation} \label{eq.L}
\sqrt{b} T  L_T^\ast \stackrel{d}{\rightarrow} \mathcal N(0,\widetilde{\theta}_0^2).
\end{equation}

  Let $I^{\ast^C}_{X,\lambda_t}=I_{X,\lambda_t}-\widehat{\mathcal F}_{\lambda_t}$ and similarly for $I^{\ast^C}_{Y,\lambda_t}$.  Verify that  
\begin{align} \label{eq.ForMean}
 {\rm E}^{\star}\langle I^{\ast^C}_{X,\lambda_t}, v_{j_1}\otimes v_{j_2}\rangle_{HS}\langle  I^{\ast^C}_{X,\lambda_t}, v_{j_1}  \otimes v_{j_2}\rangle_{HS} 
 &=  {\rm E}^{\star}\langle I^{\ast^C}_{X,\lambda_t}(v_{j_2}),v_{j_1}\rangle \langle I^{\ast^C}_{X,\lambda_t}(v_{j_2}),v_{j_1}\rangle \nonumber \\
 &= \langle {\rm E}^{\star} I^{\ast^C}_{X,\lambda_t}(v_{j_2})\otimes \overline{I}^{\ast^C}_{X,\lambda_t}(v_{j_1}), 
 v_{j_1}\otimes v_{j_2}\rangle_{HS} \nonumber \\
 & = \langle \widehat{\mathcal F}_{\lambda_t}(v_{j_1})\otimes \widehat{\mathcal F}_{\lambda_t}(v_{j_2}), 
 v_{j_1}\otimes v_{j_2}\rangle_{HS} \nonumber \\
 &= \langle \widehat{\mathcal F}_{\lambda_t}(v_{j_1}) , v_{j_1} \rangle  \langle v_{j_2}, \widehat{\mathcal F}_{\lambda_t}(v_{j_2}) \rangle.
 \end{align}
 Furthermore, 
 \begin{align} \label{eq.ForVar}
 {\rm cov}^\ast( D^*_{t}(j_1, j_2), &D^*_{t}(r_1, r_2) )  ={\rm E}^{\star} ( D^*_{t}(j_1, j_2) \overline{D}^*_{t}(r_1, r_2)) \nonumber \\
  &= {\rm E}^{\star} \langle I^{\ast^C}_{X,\lambda_t}, v_{j_1}\otimes v_{j_2}\rangle_{HS}\langle \overline{I}^{\ast^C}_{X,\lambda_t}, v_{r_1}\otimes v_{r_2}\rangle_{HS}   + {\rm E}^{\star}\langle I^{\ast^C}_{Y,\lambda_t}, v_{j_1}\otimes v_{j_2}\rangle_{HS}\langle \overline{I}^{\ast^C}_{Y,\lambda_t}, v_{r_1}\otimes v_{r_2}\rangle_{HS} \nonumber\\
 & = 2 \langle \widehat{\mathcal F}_{\lambda_t}(v_{r_2}) , v_{j_1} \rangle  \langle v_{j_2}, \widehat{\mathcal F}_{\lambda_t}(v_{r_1}) \rangle = 2 \langle \widehat{\mathcal F}_{\lambda_t}(v_{r_2}) \otimes\widehat{\mathcal F}_{\lambda_t}(v_{r_1}) , v_{j_1} \otimes v_{j_2}  \rangle_{HS},
 \end{align}
 where the last two equalities follow using the derivations in  (\ref{eq.ForMean}).

Consider  first (\ref{eq.M}). Using  (\ref{eq.ForMean}), we get 
\begin{align} 
{\rm E}^\ast(\sqrt{b} T M^\ast_T)  &=  \frac{2\pi}{T^2 b^{3/2}} \sum_{j_1,j_2=1}^\infty \sum_{t=-N}^N\sum_{l=-N}^N W^2\Big(\frac{\lambda_l-\lambda_{t}}{b} \Big){\rm E}^{\star}\langle I^{\ast^C}_{X,\lambda_t}-I^{\ast^C}_{Y,\lambda_t}, v_{j_1}\otimes v_{j_2}\rangle^2_{HS} \nonumber \\
&=  \frac{2\pi}{T^2 b^{3/2}} \sum_{j_1,j_2=1}^\infty \sum_{t=-N}^N\sum_{l=-N}^N W^2\Big(\frac{\lambda_l-\lambda_{t}}{b} \Big) \big\{
{\rm E}^{\star}\langle I^{\ast^C}_{X,\lambda_t}, v_{j_1}\otimes v_{j_2}\rangle_{HS}\langle I^{\ast^C}_{X,\lambda_t}, v_{j_1}\otimes v_{j_2}\rangle_{HS} \nonumber \\
& \ \ \ \ \ \ \ \ \ \ + {\rm E}^{\star}\langle I^{\ast^C}_{Y,\lambda_t}, v_{j_1}\otimes v_{j_2}\rangle_{HS}\langle I^{\ast^C}_{Y,\lambda_t}, v_{j_1}\otimes v_{j_2}\rangle_{HS}\big\} \nonumber \\
& = \frac{4\pi}{T^2 b^{3/2}} \sum_{j_1,j_2=1}^\infty \sum_{t=-N}^N\sum_{l=-N}^N W^2\Big(\frac{\lambda_l-\lambda_{t}}{b} \Big)\langle \widehat{\mathcal F}_{X,\lambda_t}(v_{j_1}) , v_{j_1} \rangle  \langle v_{j_2}, \widehat{\mathcal F}_{X,\lambda_t}(v_{j_2}) \rangle\nonumber \\
%& =\frac{4\pi}{T^2 b^{3/2}} \sum_{j_1,j_2=1}^\infty \sum_{t=-N}^N\sum_{l=-N}^N W^2\Big(\frac{\lambda_l-\lambda_{t}}{b} \Big)
%\langle (\widehat{\mathcal F}_{\lambda_t}\otimes \widehat{\mathcal F}_{\lambda_t})(v_{j_1}\otimes v_{j_2}) , v_{j_1}\otimes v_{j_2}\rangle_{HS} \nonumber \\
& = \frac{4\pi}{T^2b^{3/2}}  \sum_{t=-N}^N\sum_{l=-N}^N W^2\big(\frac{\lambda_l-\lambda_{t}}{b} \big) \big(\mathrm{trace}( \widehat{\mathcal F}_{\lambda_t}) \big)^2 \nonumber \\
& = \frac{4\pi}{T^2b^{3/2}}  \sum_{t=-N}^N\sum_{l=-N}^N W^2\big(\frac{\lambda_l-\lambda_{t}}{b} \big) \big(\mathrm{trace}( {\mathcal F}_{\lambda_t} )\big)^2 + o_P(1). \nonumber
\end{align}
and, therefore, 
\begin{equation} \label{eq.MeanBoot}
 b^{1/2}{\rm E}^\ast(\sqrt{b} T M^\ast_T)  = \frac{4\pi}{T^2b}  \sum_{t=-N}^N\sum_{l=-N}^N W^2\big(\frac{\lambda_l-\lambda_{t}}{b} \big) \big(\mathrm{trace}( F_{\lambda_t} )\big)^2 +o_P(1) \stackrel{P}{\rightarrow} \widetilde{\mu}_0.  
 \end{equation}
Furthermore,
\begin{align*}
{\rm var}^\ast(\sqrt{b} T M^\ast_T) & = \frac{4\pi^2}{T^4b^3}\sum_{j_1,j_2=1}^\infty\sum_{r_1,r_2=1}^\infty \sum_{t_1,t_2=-N}^N \sum_{l_1,l_2=-N}^N W^2\Big(\frac{\lambda_{l_1}-\lambda_{t_1}}{b} \Big)
W^2\Big(\frac{\lambda_{l_2}-\lambda_{t_2}}{b} \Big)\\
& \ \ \  \ \times {\rm cov}^\ast(D^\ast_{{t_1}}(j_1,j_2), D^\ast_{{t_2}}(r_1,r_2) )
\end{align*} 
which due to the independence of $D^\ast_{{t_1}}(j_1,j_2) $ and $D^\ast_{{t_2}}(j_1,j_2) $ for $ |\lambda_{t_1}| \neq |\lambda_{t_2}|$, is  reduced to four terms 
with a typical one given by  
$$
\frac{4\pi^2}{T^4b^3}  \sum_{t=1}^N \sum_{l_1,l_2=-N}^N W^2\Big(\frac{\lambda_{l_1}-\lambda_{t}}{b} \Big)
W^2\Big(\frac{\lambda_{l_2}-\lambda_{t}}{b} \Big)  \ \times \sum_{j_1,j_2=1}^\infty\sum_{r_1,r_2=1}^\infty  {\rm cov}^\ast(D^\ast_{{t}}(j_1,j_2), D^\ast_{{t}}(r_1,r_2) )
$$
and which is easily seen to  be of order  $ O_P((Tb)^{-1})$.
Similar arguments  applied to the  other three terms show that they also are asymptotically negligible  from which we get that   ${\rm var}^\ast(\sqrt{b} T M^\ast_T) \stackrel{P}{\rightarrow} 0$.
In view of (\ref{eq.MeanBoot}) this implies  that   $ \sqrt{b} T M_T^\ast \,-\, b^{-1/2}\,\widetilde{\mu}_0 \stackrel{P}{\rightarrow} 0$. 

Consider next  (\ref{eq.L}).  Notice  that
\begin{align*}
&{\rm var}^\ast(\sqrt{b}T L^\ast_T) = \frac{4\pi^2}{T^4b^3}\sum_{j_1,j_2=1}^\infty \sum_{r_1,r_2=1}^\infty \sum_{\substack{t_1,t_2=-N \\ t_1\neq t_2}}^N\sum_{\substack{t_3,t_4=-N \\ t_3\neq t_4}}^N\sum_{l_1,l_2 =-N}^N W\Big(\frac{\lambda_{l_1}-\lambda_{t_1}}{b}\Big)W\Big(\frac{\lambda_{l_1}-\lambda_{t_2}}{b}\Big)\\
& \quad  \times W\Big(\frac{\lambda_{l_2}-\lambda_{t_3}}{b}\Big)W\Big(\frac{\lambda_{l_2}-\lambda_{t_4}}{b}\Big) \Big\{{\rm E}^\ast\big(D^\ast_{t_1}(j_1,j_2) \overline{D^\ast_{t_3}}(r_1,r_2)\big) 
{\rm E}^\ast\big(D^\ast_{t_2}(j_1,j_2)\overline{D^\ast_{t_4}}(r_1,r_2)\big)\\
& + {\rm E}^\ast\big(D^\ast_{t_1}(j_1,j_2)\overline{D^\ast_{t_4}}(r_1,r_2)\big) 
{\rm E}^\ast\big(D^\ast_{t_2}(j_1,j_2)\overline{D^\ast_{t_3}}(r_1,r_2)\big) + {\rm cum}^{\star}\big(D^\ast_{t_1}(j_1,j_2), \overline{D^\ast_{t_2}}(j_1,j_2),D^\ast_{t_3}(r_1,r_2),
\overline{D^\ast_{t_4}}(r_1,r_2) \big)\Big\} \\
& = V_{1,T}^\ast + V_{2,T}^\ast + V_{3,T}^\ast,
\end{align*}
with an obvious notation for $V_{i,T}^\ast$, $i \in \{1,2,3\}$.  Since   $ {\rm E}^\ast\big(D^\ast_{t}(j_1,j_2)\overline{D^\ast_{s}}(r_1,r_2)\big)=0$ for $ |t|\neq |s|$  we get 
using (\ref{eq.ForVar})  and 
%\[
$\sum_{j_1,j_2=1}^\infty  \langle \widehat{\mathcal F}_{\lambda_{t_1}}(v_{r_2}) \otimes \widehat{\mathcal F}_{\lambda_{t_1}}(v_{r_1}),  v_{j_1}\otimes v_{j_2}   \rangle_{HS}(v_{j_1}\otimes v_{j_2} )= 
  \widehat{\mathcal F}_{\lambda_{t_1}}(v_{r_2}) \otimes \widehat{\mathcal F}_{\lambda_{t_1}}(v_{r_1})$,
  %   \]
that 
%and  that
%\[  E^\ast\big(D^\ast_{t}(j_1,j_2)\overline{D^\ast_{t}}(r_1,r_2)\big)= 2 \langle (\widehat{\mathcal F}_{\lambda_t} \otimes \widehat{\mathcal F}_{\lambda_t} ) (v_{j_1}\otimes v_{j_2}) , v_{r_1}\otimes %v_{r_2}  \rangle_{HS}.\] 
%Therefore,  
\begin{align*}
& V_{1,T}^\ast =\frac{16\pi^2}{T^4b^3}\sum_{j_1,j_2=1}^\infty \sum_{r_1,r_2=1}^\infty \sum_{t_1,t_2=-N}^N \sum_{l_1,l_2=-N}^N 
W\Big(\frac{\lambda_{l_1}-\lambda_{t_1}}{b}\Big)  W\Big(\frac{\lambda_{l_2}-\lambda_{t_1}}{b}\Big)W\Big(\frac{\lambda_{l_2}-\lambda_{t_2}}{b}\Big)W\Big(\frac{\lambda_{l_1}-\lambda_{t_2}}{b}\Big)\\
&  \quad \times \langle \widehat{\mathcal F}_{\lambda_{t_1}}(v_{r_2}) \otimes \widehat{\mathcal F}_{\lambda_{t_1}}(v_{r_1}),  v_{j_1}\otimes v_{j_2}   \rangle_{HS}  \langle \widehat{\mathcal F}_{\lambda_{t_2}}(v_{r_2}) \otimes \widehat{\mathcal F}_{\lambda_{t_2}}  (v_{r_1}),  v_{j_1}\otimes v_{j_2} \rangle_{HS}\\
%&=\frac{16\pi^2}{T^4b^3}\sum_{r_1,r_2=1}^\infty \sum_{t_1,t_2=-N}^N \sum_{l_1,l_2=-N}^N 
%W\Big(\frac{\lambda_{l_1}-\lambda_{t_1}}{b}\Big)\\
%&\ \ \ \ \times W\Big(\frac{\lambda_{l_2}-\lambda_{t_1}}{b}\Big)W\Big(\frac{\lambda_{l_2}-\lambda_{t_2}}{b}\Big)W\Big(\frac{\lambda_{l_1}-\lambda_{t_2}}{b}\Big)\\
%& \ \ \ \ \times \langle \widehat{\mathcal F}_{\lambda_{t_2}}(v_{r_2})\otimes  \widehat{\mathcal F}_{\lambda_{t_2}}(v_{r_1}), 
%    \widehat{\mathcal F}_{\lambda_{t_1}}(v_{r_2}) \otimes \widehat{\mathcal F}_{\lambda_{t_1}}  (v_{r_1})\rangle_{HS} \\
&=\frac{16\pi^2}{T^4b^3}\sum_{r_1,r_2=1}^\infty \sum_{t_1,t_2=-N}^N \sum_{l_1,l_2=-N}^N 
W\Big(\frac{\lambda_{l_1}-\lambda_{t_1}}{b}\Big)  W\Big(\frac{\lambda_{l_2}-\lambda_{t_1}}{b}\Big)W\Big(\frac{\lambda_{l_2}-\lambda_{t_2}}{b}\Big)W\Big(\frac{\lambda_{l_1}-\lambda_{t_2}}{b}\Big)\\
& \ \ \ \ \times \langle \widehat{\mathcal F}_{\lambda_{t_2}}(v_{r_2}) , \widehat{\mathcal F}_{\lambda_{t_1}}(v_{r_2}) \rangle \langle  \widehat{\mathcal F}_{\lambda_{t_1}}(v_{r_1}),
 \widehat{\mathcal F}_{\lambda_{t_2}}  (v_{r_1})\rangle\\
%& \ \ \ \ \times \langle (\widehat{\mathcal F}_{\lambda_{t_2}} \otimes \widehat{\mathcal F}_{\lambda_{t_2}} ) (v_{j_1}\otimes v_{j_2}) , v_{r_1}\otimes v_{r_2}  \rangle_{HS}\\
& = \frac{16\pi^2}{T^4b^3} \sum_{t_1,t_2=-N}^N \sum_{l_1,l_2=-N}^N 
W\Big(\frac{\lambda_{l_1}-\lambda_{t_1}}{b}\Big) W\Big(\frac{\lambda_{l_2}-\lambda_{t_1}}{b}\Big)W\Big(\frac{\lambda_{l_2}-\lambda_{t_2}}{b}\Big)W\Big(\frac{\lambda_{l_1}-\lambda_{t_2}}{b}\Big)
\langle \widehat{\mathcal F}_{\lambda_{t_1}},  \widehat{\mathcal F}_{\lambda_{t_2}}\rangle_{HS}^2 \\
& = \frac{4}{T^2b^3} \sum_{t_1, t_2=-N}^N \Big(\frac{2\pi}{T}\sum_{l=-N}^N W\Big(\frac{\lambda_{l}-\lambda_{t_1}}{b}\Big)W\Big(\frac{\lambda_{l}-\lambda_{t_2}}{b}\Big) \Big)^2 \langle \widehat{\mathcal F}_{\lambda_{t_1}},  \widehat{\mathcal F}_{\lambda_{t_2}}\rangle_{HS}^2 \\
& \rightarrow \frac{2}{\pi^2} \int_{-2\pi}^{2\pi}\Big( \int_{-\pi}^{\pi} W(u) W(u-x) du\Big)^2 dx \int_{-\pi}^\pi \|{\mathcal F}_\lambda\|^4 d\lambda,
\end{align*}
where the last convergence follows by the same arguments as in proving assertion (i) appearing in the proof of Lemma 3 in the Supplementary Material. 

Along the same lines, the same  expression is obtained  for the probability limit of $ V_{2,T}^\ast$,  while under the assumptions made,  $ V_{3,T}^\ast \rightarrow 0$ in probability. To see why the last statement is true,  use the  notation 
\[ w(i,j,k,l) = W\Big(\frac{\lambda_{l_i}-\lambda_{k}}{b}\Big)W\Big(\frac{\lambda_{l_j}-\lambda_{k}}{b}\Big)
W\Big(\frac{\lambda_{l_i}-\lambda_{l}}{b}\Big)W\Big(\frac{\lambda_{l_j}-\lambda_{l}}{b}\Big) ,\]
and observe that   $ D^\ast_{-t}(j_1,j_2) = D^\ast_t(j_1, j_2)$. By the independence of the random variables 
$D^\ast_{t}(j_1,j_2) $ and $D^\ast_{s}(j_1,j_2) $ for frequencies $ |t|\neq |s|$,  we get 
that 
\begin{align*}
V_{3,T}^\ast  &= 
\frac{1}{T^4b^3}\sum_{j_1,j_2=1}^\infty \sum_{r_1,r_2=1}^\infty \sum_{\substack{t_1,t_2=-N \\ t_1\neq t_2}}^N\sum_{l_1,l_2 =-N}^N 
w( l_1,l_2,t_1,t_2) {\rm cum}^\ast\Big(D^\ast_{t_1}(j_1,j_2), \overline{D}^\ast_{t_1}(r_1,r_2) , \overline{D}^\ast_{t_2}(j_1,j_2), D^\ast_{t_2}(r_1,r_2\Big) \\
%\end{align*}
%\\
%\begin{align*}
& = \frac{1}{T^4b^3}  \sum_{\substack{t_1, t_2=1 \\t_1\neq t_2}}^N\sum_{l_1=-N}^N\sum_{l_2=-N}^N\Big\{w( l_1,l_2,-t_1,-t_2) 
{\rm cum}^\ast\Big(D^\ast_{t_1}(j_1,j_2), \overline{D}^\ast_{t_1}(r_1,r_2) , \overline{D}^\ast_{t_2}(j_1,j_2), D^\ast_{t_2}(r_1,r_2\Big)\\
& \ \ \ \ \ \ + w( l_1,l_2,-t_1,t_2){\rm cum}^\ast\Big(D^\ast_{t_1}(j_1,j_2), \overline{D}^\ast_{t_1}(r_1,r_2) , \overline{D}^\ast_{t_2}(j_1,j_2), D^\ast_{t_2}(r_1,r_2\Big)\\
& \ \ \ \ \ \ + w( l_1,l_2,t_1,-t_2) {\rm cum}^\ast\Big(D^\ast_{t_1}(j_1,j_2), \overline{D}^\ast_{t_1}(r_1,r_2) , \overline{D}^\ast_{t_2}(j_1,j_2), D^\ast_{t_2}(r_1,r_2\Big)\\
& \ \ \ \ \ \ + w( l_1,l_2,t_1,t_2) {\rm cum}^\ast\Big(D^\ast_{t_1}(j_1,j_2), \overline{D}^\ast_{t_1}(r_1,r_2) , \overline{D}^\ast_{t_2}(j_1,j_2), D^\ast_{t_2}(r_1,r_2\Big)\Big\}
\end{align*}
which vanishes due to  the independence of the bootstrap finite Fourier transforms  and consequently of the random variables 
$ D^\ast_{t_1}(\cdot)$ and $D^\ast_{t_2}(\cdot) $ for $ 1\leq t_1 \neq t_2 \leq N$. 

We next show that $ \sqrt{b} T L^\ast_T \stackrel{D}{\rightarrow} \mathcal{N}(0, \widetilde{\theta}_0)$. Toward this we write 
%\[
$ \sqrt{b} T L^\ast_{T} = $ $  \sum_{j_1,j_2=1}^\infty \sum_{1\leq t_1 < t_2 \leq N} $ $  H^\ast_{t_1,t_2}(j_1,j_2),
$
%\]
where
\begin{align}\label{eq.Hstar}
H_{t_1,t_2}^\ast(j_1,j_2)&= 2\Big\{ h^\ast_{t_1,t_2}(j_1,j_2)+ h^\ast_{-t_1,t_2}(j_1,j_2) +h^\ast_{t_1,-t_2}(j_1,j_2) + h^\ast_{-t_1,-t_2}(j_1,j_2) \Big\}
\end{align}
and
\[ h_{t,s}^\ast(j,r) = \frac{2\pi}{b^{3/2} T^2}\sum_{l=-N}^N 
W\Big(\frac{\lambda_{l}-\lambda_{t}}{b}\Big)W\Big(\frac{\lambda_{l}-\lambda_{s}}{b}\Big)D_{t}^\ast(j,r) D^\ast_s(j,r).\]
Let $ \sqrt{b}T L_{T,K}^\ast = \sum_{j_1,j_2=1}^K \sum_{1\leq t_1 < t_2 \leq N} H_{t_1,t_2}^\ast(j_1,j_2)$ and 
\begin{align*}
&\widetilde{\theta}^2_{0,K}  =  \frac{4}{\pi^2} \int_{-2\pi}^{2\pi}\Big( \int_{-\pi}^{\pi} W(u) W(u-x) du\Big)^2 dx \sum_{j_1,j_2,r_1,r_2=1}^K \int_{-\pi}^\pi \langle  v_{j_1}\otimes v_{j_2}, {\mathcal F}_\lambda\rangle^2_{HS}  \langle  
v_{r_1}\otimes v_{r_2}, {\mathcal F}_\lambda\rangle^2_{HS} d\lambda \,.
\end{align*}
Then, to  establish the desired weak convergence it suffices to  prove that  \\
%\begin{enumerate}
%\item[(i)] \ 
(i) $ \sqrt{b} T L^\ast_{T,K} \stackrel{D}{\rightarrow} \mathcal{N}(0, \widetilde{\theta}_{0,K}^2)$  as $n\rightarrow \infty$ for every $K\in\N$,\\
%\item[(ii)] \ 
(ii) $ \widetilde{\theta}_{0,K}^2 \rightarrow \widetilde{\theta}_{0}^2$ as $ K\rightarrow \infty$,\\
%\item[(iii)] \ 
(iii) For every $\epsilon >0$, $ \lim_{K\rightarrow \infty} \limsup_n \Pr\big(\big| \sqrt{b}T L^\ast_{T,K} -\sqrt{b}T L_T^\ast  \big| > \epsilon \big) = 0$.\\
%\end{enumerate}
Consider (i). Observe that $  \sqrt{b} T L^\ast_{T,K}$ is a quadratic form in  the independent random variables $ D_{t}(i,j) $ and $ D_s(i,j)$, $t\neq s$. We can, therefore,  use Theorem 2.1 of de Jong \cite{dJ-1987} to establish  the weak convergence  
(i). For this we need to show that 
\begin{itemize}
\item[(a)] \   $\sigma^{-2}(T) \max_{1\leq i \leq N} \sum_{1\leq j \leq N} \sigma^2_{i,j} \rightarrow 0$, 
\item[(b)] \  ${\rm E}^{\ast}\Big( \sum_{j_1,j_2=1}^K \sum_{1\leq t_1<t_2 \leq N} H^\ast_{t_1,t_2}(j_1,j_2)\Big)^4/\sigma^4(T) \rightarrow 0$,
\end{itemize}
in probability as $ T \rightarrow \infty$, 
where $ \sigma^2(T)=\sum_{1\leq t_1 < t_2 \leq N}\sigma^2_{t_1,t_2}$ and 
\[ \sigma^2_{t_1,t_2}=\sum_{j_1,j_2,r_1,r_2=1}^K 
{\rm cov}^\ast(H_{t_1,t_2}^\ast(j_1,j_2),H_{t_1,t_2}^\ast(r_1,r_2)).\]
Evaluating  $ \sigma_{t_1,t_2}^2 = {\rm E}^\ast(\sum_{j_1,j_2=1}^K H^\ast_{t_1,t_2}(j_1,j_2))^2$  for $ 1 \leq t_1 < t_2 \leq N$, using   (\ref{eq.Hstar}), 
yields  the expression 
\[ 4\sum_{j_1,j_2,r_1,r_2=1}^K \sum_{m_1\in\{-t_1,t_1\}}\sum_{s_1\in\{-t_2,t_2\}}\sum_{m_2\in\{-t_1,t_1\}}\sum_{s_2\in\{-t_2,t_2\}} {\rm cov}^\ast(h^\ast_{m_1,s_1}(j_1,j_2),h^\ast_{m_2,s_2}(r_1,r_2)).\]
 Taking into account  
the independence of the random variables 
involved,  ($ t_1\neq t_2$),  the covariance terms  in the above sum  are very similar with a typical one given, for instance for $m_1=t_1,s_1=t_2,m_2=-t_1, s_2=-t_2 $, by 
\begin{align*}
&\frac{1}{T^4b^3}\sum_{l_1}\sum_{l_2}  W\Big(\frac{\lambda_{l_1}-\lambda_{t_1}}{b}\Big)W\Big(\frac{\lambda_{l_1}-\lambda_{t_2}}{b}\Big)W\Big( \frac{\lambda_{l_2}+\lambda_{t_1}}{b}\Big) W\Big(\frac{\lambda_{l_2}+\lambda_{t_2}}{b} \Big) \\
& \ \ \times \langle \widehat{\mathcal F}_{\lambda_{t_1}}(v_{r_2}) \otimes \widehat{\mathcal F}_{-\lambda_{t_1}}  (v_{r_1}), v_{j_1}  \otimes  v_{j_2}  \rangle_{HS} \langle \widehat{\mathcal F}_{\lambda_{t_2}}(v_{r_2}) \otimes \widehat{\mathcal F}_{-\lambda_{t_2}}  (v_{r_1}), v_{j_1}  \otimes  v_{j_2}  \rangle_{HS}\\
%&
%\times 
%\widehat{f}_{\lambda_{t_1}}(\sigma_{l_1},\tau_{l_1})\overline{ \widehat{f}}_{\lambda_{-t_1}}(\sigma_{l_2},\tau_{l_2})  \widehat{f}_{\lambda_{t_2}}(\sigma_{l_1},\tau_{l_1})  \widehat{f}_{\lambda_{-t_2}}(\sigma_{l_2},\tau_{l_2}) \\
& = \frac{1}{4\pi^2T^2b^3} \Big(\frac{2\pi}{T}\sum_{l=-N}^N W\Big(\frac{\lambda_{l}-\lambda_{t_1}}{b}\Big)W\Big(\frac{\lambda_{l}-\lambda_{t_2}}{b}\Big) \Big) \Big(\frac{2\pi}{T}\sum_{l=-N}^N W\Big(\frac{\lambda_{l}+\lambda_{t_1}}{b}\Big)W\Big(\frac{\lambda_{l}+\lambda_{t_2}}{b}\Big) \Big) \\
& \ \ \times \langle \widehat{\mathcal F}_{\lambda_{t_1}}(v_{r_2}) \otimes \widehat{\mathcal F}_{-\lambda_{t_1}}  (v_{r_1}), v_{j_1}  \otimes  v_{j_2}  \rangle_{HS}\langle \widehat{\mathcal F}_{\lambda_{t_2}}(v_{r_2}) \otimes \widehat{\mathcal F}_{-\lambda_{t_2}}  (v_{r_1}), v_{j_1}  \otimes  v_{j_2}  \rangle_{HS} = O_P(T^{-2}b^{-1}), 
\end{align*}
where the $O_P(T^{-2}b^{-1})$ term is uniform in $ t_1$ and $t_2$ because   
\[ | \langle \widehat{\mathcal F}_{\lambda_{t_1}}(v_{r_2}) \otimes \widehat{\mathcal F}_{-\lambda_{t_1}}  (v_{r_1}), v_{j_1}  \otimes  v_{j_2}  \rangle_{HS}| \leq \|\widehat{\mathcal F}_{\lambda_{t_1}}\|_{HS} \|\widehat{\mathcal F}_{-\lambda_{t_1}}\|_{HS} =O_P(1),  \]
   uniformly in $t_1$, $t_2$, and 
\begin{align*}
 &\frac{2\pi}{T}\sum_{l=-N}^N   W\Big(\frac{\lambda_{l}-\lambda_{t_1}}{b}\Big)W\Big(\frac{\lambda_{l}-\lambda_{t_2}}{b}\Big)= \int 
 W\Big(\frac{\lambda-\lambda_{t_1}}{b}\Big)W\Big(\frac{\lambda-\lambda_{t_2}}{b}\Big) +O(T^{-1})\\
 & = b  \int W\Big(u- \frac{\lambda_{t_1}}{b}\Big)W\Big(u-\frac{\lambda_{t_2}}{b}\Big)du +O(T^{-1}) = b  \int W\Big(x- \frac{\lambda_{t_1}-\lambda_{t_2}}{b}\Big)W\Big(x\Big)dx +O(T^{-1}) = O(b),
 \end{align*}
 uniformly in $t_1$, $t_2$.  
 Taking into account that  $ 0< \sigma^2(T)={\rm E}^\ast(  \sum_{j_1,j_2=1}^K  \sum_{1\leq t_1 < t_2 \leq N}  H^\ast_{t_1,t_2}(j_1,j_2) )^2=O_P(1)$, which follows from the calculations of  ${\rm var}^\ast(\sqrt{b} T L^\ast_T) $,  we get that 
\[ \frac{1}{\sigma^2(T)} \max_{1\leq t_1 \leq  N}\sum_{1\leq t_2\leq N} \sigma^2_{t_1,t_2} =  O_P(T^{-1}b^{-1}) \rightarrow 0,\]
as $ T \rightarrow \infty$, which establishes  (a). 

Consider Condition (b). From  (\ref{eq.Hstar}),  the fourth moment  of \\ $ \sum_{j_1,j_2=1}^K  \sum_{1\leq t_1 < t_2 \leq N}  H^\ast_{t_1,t_2}(j_1,j_2)$  equals 
\begin{align*}
16& \sum_{j_1, \ldots,j_8=1}^K \sum_{1\leq t_1 <t_2 \leq N}  \sum_{1\leq t_3 <t_4 \leq N}\sum_{1\leq t_5 <t_6\leq N}\sum_{1\leq t_7  < t_8 \leq N}
\sum_{r_1 \in\{ -t_1,t_1\} \atop  r_2\in\{-t_2,t_2\}} \sum_{k_1 \in\{ -t_3,t_3\} \atop k_2\in\{-t_4,t_4\}} \\
&\times \sum_{n_1 \in\{ -t_5,t_5\} \atop  n_2\in\{-t_6,t_6\}} 
\sum_{v_1 \in\{ -t_7,t_7\} \atop v_2\in\{-t_8,t_8\}}{\rm E}^\ast\Big( h^\ast_{r_1,r_2}(j_1,j_2)h^\ast_{k_1,k_2}(j_3,j_4)h^\ast_{n_1,n_2}(j_5,j_6)h^\ast_{v_1,v_2}(j_7,j_8)\Big),
\end{align*}
where  only for the following four cases  the  expectation term is different from zero:  1) $(r_1,r_2)=(k_1,k_2) \neq (n_1,n_2)=(v_1,v_2)$,  
2) $(r_1,r_2)=(n_1,n_2) \neq (k_1,k_2)=(v_1,v_2)$, 3) $(r_1,r_2)=(v_1,v_2) \neq (k_1,k_2)=(n_1,n_2)$ and 4) $(r_1,r_2)=(k_1,k_2) = (n_1,n_2)=(v_1,v_2)$ and where the  notation $ (i,j)=(l,k)$  means $ i=l$ and $j=k$. Straightforward calculations  show that case 4) vanishes asymptotically while cases 1), 2) and 3) converge  to the same limit as $ \sigma^4(T)$ converges, from which we conclude  assertion (b).

Condition (ii) follows immediately from the fact that, as $ K \rightarrow \infty$,  
\[ \sum_{j_1,j_2=1}^K \langle v_{j_1}\otimes v_{j_2}, {\mathcal F}_\lambda \rangle^2_{HS} \rightarrow
\sum_{j_1,j_2=1}^ \infty \langle v_{j_1}\otimes v_{j_2}, {\mathcal F}_\lambda \rangle^2_{HS} = \|{\mathcal F}_\lambda\|^2_{HS}.\]

Finally to establish the validity of condition (iii) notice that  
\begin{align*}
\sqrt{b} T ( & L_T^\ast -L^\ast_{T,K})  = \sqrt{b}T\Big( \sum_{j_1=1}^K\sum_{j_2=K+1}^\infty \sum_{1\leq t_1<t_2\leq N} H^\ast_{t_1,t_2}(j_1,j_2) \\
& \ \  + 
\sum_{j_1=K+1}^\infty\sum_{j_2=1}^K \sum_{1\leq t_1<t_2\leq N} H^\ast_{t_1,t_2}(j_1,j_2) +  \sum_{j_1=K+1}^\infty\sum_{j_2=K+1}^\infty \sum_{1\leq t_1<t_2\leq N} H^\ast_{t_1,t_2}(j_1,j_2)\Big) = \sum_{r=1}^3 Q_{r,T}^\ast,
\end{align*}
with an obvious notation for $ Q^\ast_{r,T}$, $r=1,2,3$.  Consider $ Q^\ast_{1,T}$. We then have
\begin{align*}
{\rm E}^\ast(Q^\ast_{1,T})^2 & = \sum_{j_1,r_1=1}^\infty\sum_{j_2,r_2=K+1}^\infty \sum_{1\leq t_1 <t_2\leq N}\sum_{1\leq s_1 <s_2\leq N}{\rm cov}^\ast(H^\ast_{t_1,t_2}(j_1,j_2), H^\ast_{s_1,s_2}(r_1,r_2) ).
\end{align*}
Now, evaluating the covariance term ${\rm cov}^\ast(H^\ast_{t_1,t_2}(j_1,j_2), H^\ast_{s_1,s_2}(r_1,r_2) ) $ as in the calculations for ${\rm var}^\ast(\sqrt{b}T L^\ast_T) $, 
using (\ref{eq.ForVar}) and the fact that $ {\mathcal F}_\lambda$ is self adjoint,  we get that
\begin{align*}    
\lim_{n\rightarrow\infty}& \sum_{1\leq t_1 <t_2\leq N}  \sum_{1\leq s_1 <s_2\leq N} {\rm cov}^\ast(H^\ast_{t_1,t_2}(j_1,j_2), H^\ast_{s_1,s_2}(r_1,r_2) )\\
& = \frac{4}{\pi^2} \int_{-2\pi}^{2\pi}\Big(\int_{-\pi}^{\pi}W(u)W(u-x)du\Big)^2 \int_{-\pi}^{\pi}\langle {\mathcal F}_\lambda(v_{r_2}),  v_{j_1}\rangle ^2
\langle {\mathcal F}_\lambda(v_{j_2}),  v_{r_1}\rangle ^2d\lambda.
 \end{align*}
 Therefore,
 \begin{align*}
\lim_{n\rightarrow\infty}  {\rm E}^\ast(Q^\ast_{1,T})^2 & =\frac{4}{\pi^2} \int_{-2\pi}^{2\pi}\Big(\int_{-\pi}^{\pi}W(u)W(u-x)du\Big)^2 \int_{-\pi}^{\pi}\Big(\sum_{j_1=1}^K       
\sum_{j_2=K+1}^\infty \langle v_{j_1}, {\mathcal F}_\lambda ( v_{j_2})\rangle ^2\Big)^2\\
& \leq \frac{4}{\pi^2} \int_{-2\pi}^{2\pi}\Big(\int_{-\pi}^{\pi}W(u)W(u-x)du\Big)^2 \int_{-\pi}^{\pi}\Big(\sum_{j_2=K+1}^\infty      
\sum_{j_1=1}^\infty \langle v_{j_1}, {\mathcal F}_\lambda ( v_{j_2})\rangle ^2\Big)^2\\
& = \frac{4}{\pi^2} \int_{-2\pi}^{2\pi}\Big(\int_{-\pi}^{\pi}W(u)W(u-x)du\Big)^2 \int_{-\pi}^{\pi}\Big(      
\sum_{j_2=K+1}^\infty \|  {\mathcal F}_\lambda ( v_{j_2})\|^2\Big)^2  \rightarrow 0,
\end{align*}
 as $ K \rightarrow \infty$ since $\lim_{K\rightarrow \infty} \sum_{j_2=K+1}^\infty \|  {\mathcal F}_\lambda ( v_{j_2})\|^2 =0$. 
%  
% \begin{align*}
% \lim_{K\rightarrow\infty} & \limsup_{n\rightarrow\infty} E^\ast(Q^\ast_{1,T})^2 =0,
% 
%   \frac{4}{\pi^2} \int_{-2\pi}^{2\pi}\Big(\int_{-\pi}^{\pi}W(u)W(u-x)du\Big)^2 \int_{-\pi}^{\pi} \sum_{j_1=1}^\infty \sum_{j_2=!langle v_{j_1} \otimes v_{j_2}, {\mathcal F}_\lambda \rangle^2_{HS}
% \langle v_{r_1} \otimes v_{r_2}, {\mathcal F}_\lambda \rangle^2_{HS} d\lambda
% \end{align*}
 %in probability. \footnote{AL: evt. sollten hier Zwischenschritte angegeben werden. SP O.k. so? AL: Warum gilt $ \|\mathcal F_{\lambda}\|_{HS}=0$? Die Null muss weg oder? Dann sollte alles so passen. } 
 By the same arguments we get  that \\ $\lim_{K\rightarrow\infty}  \limsup_{n\rightarrow\infty} {\rm E}^\ast(Q^\ast_{2,T})^2 =0$ and  
 $  \lim_{K\rightarrow\infty}  \limsup_{n\rightarrow\infty} {\rm E}^\ast(Q^\ast_{3,T})^2=0$, in probability.
 Condition (iii) follows then 
using the bound   $\sqrt{b} T {\rm E}^\ast( L_T^\ast -L^\ast_{T,K})^2  \leq C \sum_{r=1}^3 {\rm E}^\ast(Q^\ast_{r,T})^2 $. 

\hspace*{0.1cm}

\noindent {\it Supplement to ``Bootstrap-Based Testing of the Equality of Spectral Density Operators for Functional  Processes''} \ The online supplement contains some useful  technical tools,  some new results on frequency domain properties of linear Hilbertian stochastic processes and 
the proofs that were omitted in this paper. 

%\section*{References}
\bibliography{}

\vspace*{0.3cm}
%

%\begin{small}
%Anne Leucht, University of Bamberg, Research Group of Statistics and Business Mathematics, Feldkirchenstra\ss e 21, D-96052 Bamberg, Germany (E-mail:  anne.leucht@uni-bamberg.de)
%
%Efstathios Paparoditis, University of Cyprus, Department of Mathematics and Statistics, P.O.Box 20537, CY-1678, Nicosia, Cyprus (stathisp@ucy.ac.cy)
%
%Daniel Rademacher, Technische Universit\"at Braunschweig, Institut f\"ur Mathematische Stochastik, Universit\"atsplatz 2, D-38106 Braunschweig, Germany (d.rademacher@tu-braunschweig.de)
%
%Theofanis Sapatinas, University of Cyprus, Department of Mathematics and Statistics, P.O.Box 20537, CY-1678, Nicosia, Cyprus (fanis@ucy.ac.cy)
%\end{small}
\end{document}